\documentclass[11pt]{article}
\usepackage{latexsym}

\usepackage{amsmath}
\usepackage{amssymb}
\usepackage{amsthm}
\usepackage{amsfonts}

\newtheorem{lemma}{Lemma}
\newtheorem{remark}{Remark}
\newtheorem{definition}{Definition}
\newtheorem{theorem}{Theorem}
\newtheorem{corollary}{Corollary}
\newtheorem{proposition}{Proposition}

\def \R{I\!\!R}

\def \N{I\!\!N}

\def \esssup{\rm{esssup}}
\def \essinf{\rm{essinf}}

\newcommand{\INTER}{\mathop{\rm {\cap}}\limits}

\begin{document}

\title{Dynamic risk measuring under model uncertainty: taking advantage of the hidden probability measure} 

\author{Jocelyne BION-NADAL
\footnote{e-mail : jocelyne.bion-nadal@cmap.polytechnique.fr}\\ \small \it  UMR 7641  CNRS CMAP
Ecole 
 Polytechnique, 91128 Palaiseau Cedex, France\\  and\\   Magali KERVAREC\footnote{ e-mail :  mkervare@univ-evry.fr}\\
\small \it 
Laboratoire analyse et probabilit\'es, Universit\'e d' Evry, Bd F. Mitterrand, 91000 Evry France}
\date{} 

 \maketitle

\begin{abstract}
We study dynamic risk measures in a very general framework enabling to model uncertainty and processes with jumps. We previously showed the existence of  a canonical equivalence class of probability measures hidden behind a given set of probability measures possibly non dominated. Taking advantage of this result, we exhibit a dual representation that completely characterizes the dynamic risk measure. We prove continuity and characterize time consistency.  Then, we prove regularity for all processes associated to time consistent convex dynamic risk measures. We also study factorization through time for sublinear risk measures. Finally we consider examples (uncertain volatility and G-expectations).

\end{abstract}


 \section{Introduction}

In a previous paper \cite{BNK} we have revealed a canonical equivalence class of probability measures hidden behind every weakly relatively compact set of probability measures, possibly non dominated. In the present paper, taking advantage of the existence of this class of probability measures, we study dynamic risk measures in a very general framework enabling model uncertainty. \\
In order to study dynamic risk measures for very general family of processes, the state space in this study will be either the space of continuous functions, or the space of c\`adl\`ag functions on $\R_+$ with values in $\R^d$. A key notion in a dynamic setting is the notion of time consistency which allows to quantify the risk in a consistent way for financial products with distinct maturity dates.  A time consistent dynamic risk measure is a family $\rho_{s,t}$ of convex conditional risk measures satisfying $\rho_{r,t}=\rho_{r,s}\circ (-\rho_{s,t})$ for every $r \leq s\leq t$. Equivalently a time consistent dynamic risk measure is a convex risk measure $\rho_{0,T}$ which can be factorized through all instants of time.  Time consistent dynamic risk measures have been studied in several papers in the setting of a filtered probability space. F. Delbaen \cite{D} has characterized time consistency for dynamic sublinear (or coherent) risk measures. For convex dynamic risk measures, time consistency has been discussed by P. Cheridito et al   \cite{CDK}, S. Kl\"oppel and M. Schweizer  \cite{KS}, J Bion-Nadal  \cite{BN02} \cite{JBN}, and B. Roorda and H. Schumacher  \cite{RS}. In particular a characterization of time consistency in terms of acceptance sets can be found in \cite{CDK}. A characterization in terms of a cocycle condition is proved in  \cite{BN02}, \cite{JBN}. This last characterization  is very tractable and allows for the construction of families of time consistent dynamic risk measures. 
S. Peng \cite{Peng}, N. El Karoui and P. Barrieu  \cite{BEK} and F. Delbaen  et al  \cite{DPR}, are focusing on the particular case of a Brownian filtration. In particular it has been proved in \cite{DPR} that every normalized time consistent dynamic risk measure on a Brownian filtration is the limit of a sequence of solutions of B S D E (Backward Stochastic Differential Equations).
 Notice that in all the previous cited papers on dynamic risk measures, a  reference probability space is fixed. This framework is sufficiently rich to include the study of models with stochastic volatility or models with jumps. However, the framework of a fixed probability space is not appropriate to study financial problems with uncertainty about the volatility, or more generally to study problems using a family of probability measures which are not dominated by a single probability measure.\\
 Some recent papers have considered this context of model uncertainty.  L. denis and C. Martini \cite{DM} have studied the problem of superhedging (equivalent to sublinear risk measuring) under uncertain volatility.  S. Peng has introduced in  \cite{P1} and \cite{P2} conditional G-expectations ${\cal E}^G_t$  defined on a subset of continuous bounded functions on the canonical space of continuous paths. More precisely $B_t$ denoting the coordinate process on $\Omega= {\cal C}([0,T],\R^d)$, ${\cal E}^G_t(\phi(B_{t_1},B_{t_2}, ...B_{t_k})$ is defined for $\phi$  Lipschitz bounded on $\R^k$ recursively on $k$ using solutions of partial differential equations (PDE). ${\cal E}^G_t$ is sublinear monotone and satisfies the time consistency property: for $s\leq t$, ${\cal E}^G_s={\cal E}^G_s \circ {\cal E}^G_t$. It is thus (up to a minus sign) a time consistent sublinear dynamic risk measure defined on a subset of ${\cal C}_b(\Omega)$.  L.Denis et al  \cite{DHP} have  proved that ${\cal E}_0$ admits the following dual representation ${\cal E}^G_0 (X)= \sup_{P \in {\cal P}}E_P(X)$. The important property of  ${\cal P}$ is that it is a weakly relatively compact set of probability measures not all dominated by a single probability measure.  
In \cite{BNK} we have proved more generally  that every regular convex  risk measure $\rho$ on ${\cal C}_b(\Omega)$ where $\Omega$ is a Polish space can be written $\rho(X)=\sup_ {P \in {\cal P}}(E_P(-X)-\alpha(P))$ for some weakly relatively compact set ${\cal P}$ of probability measures on $(\Omega,{\cal B}(\Omega))$.
The  very recent papers of H.M.Soner et al \cite{STZ} and M. Nutz \cite{Nutz}  construct examples of time consistent dynamic risk measures in this uncertainty framework, without making use of PDE.  In  these papers, as in the case of G-expectations, $\Omega={\cal C}([0,T],\R^d)$,  and  ${\cal P}$ is always a subset of the laws $P^{\theta}=P^0 \circ (X^{\theta})^{-1}$, where $P^0$ is the Wiener measure and  $X^{\theta}_t=\int _0^t \theta(s)dB_s$ , $\theta$ taking values in the positive definite matrices. ${\cal E}_0(X)=\sup_{P \in {\cal P}} E_P(X)$. For $X$ belonging to ${\cal U}{\cal C}_b(\Omega)$, the set of uniformly continuous bounded function on $\Omega$, ${\cal E}_t(X)(\omega)$ is defined pointwise. In \cite{STZ} time consistent dynamic risk measure are constructed on ${\cal U}{\cal C}_b(\Omega)$ with values in ${\cal U}{\cal C}_b(\Omega_t)$ where $\Omega_t={\cal C}([0,t],\R^d)$ making use of solutions of B S D E. Time consistent sublinear dynamic risk measures are constructed in \cite{Nutz},   making use of regular conditional probability distributions. \\
The goal  of the present paper is to conduct a systematic study of dynamic risk measures in a very general framework of uncertainty represented by a weakly relatively compact set of probability measures on a Polish space $\Omega$. This study  is in the line of the previous study in  \cite{BN02} of dynamic risk measures on a fixed probability space.  In order to deal with general families of processes, $\Omega$ will be either the space of continuous paths, or the space of c\`adl\`ag paths.
In our previous paper \cite{BNK} we have proved that  every convex regular risk measure on ${\cal C}_b(\Omega)$ admits  a dual representation with a  weakly relatively compact set of probability measures ${\cal P}$ and extends thus uniquely to  $L^1(c)$, the Banach space obtained from completion and separation of ${\cal C}_b(\Omega)$  with respect to the  semi-norm $c$, $c(f)=\sup_{P \in{\cal P}}E_P(|f|)$. A key result of \cite{BNK} is the existence of an equivalence class of probability measures, such that if $P$ belongs to this class, every element $X$ in $L^1(c)$ is completely determined by its value $P$ a.s.  From a dynamic point of view, it would be  too restrictive to study dynamic risk measures taking values in the space of continuous functions ${\cal C}_b(\Omega_t)$ where $\Omega_t$ denotes either the space of continuous paths ${\cal C}([0,t],\R^d)$ or  the space of c\`adl\`ag paths ${\cal D}([0,t[,\R^d)$. Therefore we introduce an axiomatic for dynamic risk measures on $L^1(c)$.\\
The paper is organized as follows. In Section 2, we introduce a uniform Skorokhod topology on compact spaces on ${\cal D}([0,t[,\R^d)$ in order to construct a continuous projection from ${\cal D}([0,\infty[,\R^d)$ onto ${\cal D}([0,t[,\R^d)$ and thus an embedding of ${\cal C}_b({\cal D}([0,t[,\R^d)$ into ${\cal C}_b({\cal D}([0,\infty[,\R^d)$. In section 3, we use the topological results of section 2 to define the increasing set of Banach spaces $\left( L^1_t(c)\right)$ and the axiomatic for dynamic convex risk measures on $L^1(c)$. In Section 4, we give a characterization of  convex dynamic risk measure admitting a  dual representation in terms of probability measures. The necessary and sufficient  condition is here a strong convexity condition (whereas in case of dynamic risk measures on $L^{\infty}$ spaces the condition was a condition of semi-continuity). Notice that in all this  study a very important point is the  fact that an element of $L^1(c)$ is competely determined by its value $P$ a.s. where $P$ is some probability measure belonging to the canonical $c$-class. 
In section 5, we prove the continuity of dynamic risk measures for the $c$-norm. 
 Section 6 is devoted to the important property of time consistency. As in the case of $L^{\infty}$ spaces (cf. \cite{BN02}),   we give a characterization of the time consistency in terms of a cocycle property and also in terms of acceptance sets. We also prove  the regularity of paths for normalized convex time consistent dynamic risk measures. This result extends to the context of model uncertainty the regularity for paths  which was proved in \cite{D} in the case of  sublinear dynamic risk measures on $L^{\infty}$ spaces, and in \cite{JBN} in the more general case of convex dynamic risk measures. Notice that in the proof of  the regularity of paths, the difficulty consists in proving  the right continuity of $s \rightarrow (\rho_{st})(X)$. 
The proofs given here follow the proofs of \cite{BN02} and \cite{JBN}, but they  are more tricky  because  $\rho_{st}$ is defined on a Banach space $L^1_t(c)$ which does not contain the ${\cal F}_t$-measurable functions. Nutz and Soner 
\cite{NS} give a similar statement for the regularity of paths in the particular case of  sublinear risk measures. Section 7 deals with different ways to address the question of the factorization through time of a convex risk measure on ${\cal C}_b(\Omega)$. We give sufficient conditions for the existence of a factorization.  The unicity of a factorization through time is studied in \cite{NS}.  
 Finally, Section 8 is dedicated to the study of examples. In light of our previous results,  we consider the example of uncertain volatility and the case of  G-expectations.

\section{The uniform Skorokhod topology on compact spaces}
\label{secSko}
In order to deal with processes with continuous paths, we consider the set ${\cal C}_0([0,\infty[,\R^d)$ of continuous functions $f$ defined on $[0,\infty[$ with values in $\R^d$ such that $f(0)=0$. The space ${\cal C}_0([0,\infty[,\R^d)$ endowed with the topology of uniform convergence on compact spaces  is a metrizable separable  complete space (cf \cite{Bil}). \\
As we want also to deal with processes with c\`adl\`ag paths we consider the Skorokhod space $D([0,\infty[,\R^d)$   i.e. the space of c\`adl\`ag functions defined on $[0,\infty[$ with values in $\R^d$.  The space $D([0,\infty[,\R^d)$ endowed with the Skorokhod topology is a metrizable separable complete space (\cite{Bil}). The subset $D_0([0,\infty[,\R^d)$ of c\`adl\`ag functions $f$ such that $f(0)=0$ is closed for the Skorokhod topology. Thus it is also a Polish space.\\

\begin{lemma}
Let $\Omega={\cal C}_0([0,\infty[,\R^d)$,  $\Omega_t={\cal C}_0([0,t],\R^d)$
and $\tilde\Omega_t={\cal C}_0([t,\infty[,\R^d)$. The spaces $\Omega$, $\Omega_t$ and $\tilde \Omega_t$ are endowed with the topology of uniform convergence on compact spaces. The product space $\Omega_t \times \tilde \Omega_t$ is endowed with the product topology.
Then the  metrizable separable complete spaces $\Omega$ and $\Omega_t \times \tilde \Omega_t$ are homeomorphic.
\label{homeo}
\end{lemma}
{\bf Proof}
Let 
\begin{eqnarray} 
\phi: \Omega & \rightarrow  & \Omega_t \times \tilde \Omega_t \nonumber\\
       x & \rightarrow & (x_{|_{[0,t]}},x_{|_{[t, \infty[}}-x(t)) \nonumber
\label{homeomorph}
\end{eqnarray}
Let 
\begin{eqnarray} 
\psi : \Omega_t \times \tilde \Omega_t & \rightarrow & \Omega \nonumber\\
        (x_1,x_2)  & \rightarrow & x_11_{[0,t[}+(x_2+x_1(t))1_{[t, \infty[} \nonumber
\label{homeomorph2}
\end{eqnarray}
It is easy to verify that $\phi$ and $\psi$ are reciprocal bijections.
As the topological spaces $ \Omega $ and $ \Omega_t \times \tilde \Omega_t $ are metrizable spaces, we verify the  continuity of $\phi$ and $\psi$  along sequences . The continuity follows from  the fact that a sequence $x_n$ tends to $x$   uniformly on compact spaces of $[0,\infty[$ if and only if  the restriction of $x_n$ to $[0,t]$ tends uniformly to the restriction of $x$ to $[0,t]$, and the restriction of  $x_n$ to $[t,\infty[$ tends uniformly on compact spaces to the restriction of $x$ to $[t,\infty[$. 
\hfill $\square $   
\begin{remark}
If we replace  in  Lemma \ref{homeo} $\;\Omega$ by  $D([0,\infty[,\R^d)$, $\Omega_t$ by  $D([0,t],\R^d)$ and $\tilde \Omega_t$ by  $D([t,\infty[,\R^d)$ endowed with the Skorokhod topologies, the map $\phi$ is no more continuous. Indeed the map $x \in D([0,\infty[,\R^d) \rightarrow  x_{|_{[0,t]}} \in D([0,t],\R^d)$ is continuous at $x$ if $x$ is continuous at $t$ but is not continuous at every point $x$.  Therefore in order to have the continuity for the projection from $D([0,\infty[,\R^d)$ onto $D([0,t[,\R^d)$, we will introduce a new topology on $D([0,t[,\R^d)$.
 \label{remsko}
\end{remark}
\begin{definition}
Let $t>0$. On $D([0,t[,\R^d)$ we define the uniform Skorokhod topology on compact spaces  as  the topology deduced from the Skorokhod topology on $D([0,\infty[,\R^d)$ using the strictly increasing homeomorphism: 
\begin{eqnarray}
\alpha_t:[0,t[ & \rightarrow & [0, \infty[ \nonumber \\
u & \rightarrow & \frac{u}{t-u} 
\label{eqalpha}
\end{eqnarray}
The uniform Skorokhod topology on compact spaces on  $D([0,t[,\R^d)$ is thus defined by the  metric
$\hat d (x,y)=d_{\infty}(x \circ(\alpha_t)^{-1}, (y \circ(\alpha_t)^{-1})$, 
where $d_{\infty}$ is the metric introduced in \cite{Bil} Chapter 3 Section 16.
\label{defsko}
\end{definition}
Denote $\Lambda_t$ the set of continuous strictly increasing maps from $[0,t[$ onto itself.
\begin{proposition}

$D([0,t[,\R^d)$ endowed with the uniform Skorokhod topology on compact spaces is  a metrizable separable complete space.\\
The sequence  $\hat d (x_n, x)$ tends to $0$ in $D([0,t[,\R^d)$ if and only if there exists elements $\gamma_n$ of $\Lambda_t$ such that 
\begin{eqnarray}
\sup_{u \in [0,t[} |\gamma_n(u)-u| \rightarrow 0 \nonumber\\
\forall m \in \N\;\;\sup_{u \in [0,t(1-\frac{1}{1+m})]} ||x_n(\gamma_n(u))-x(u)|| \rightarrow 0 
\label{eqSkor}
\end{eqnarray}
\label{propsko}
\end{proposition}
{\bf Proof} 
 The first part of the proposition follows  easily from Theorem 16.3 of \cite{Bil}. \\
Assume that $\hat d (x_n, x) \rightarrow 0$. From Theorem 16.1 of \cite{Bil}, there is a sequence $\lambda_n \in \Lambda_{\infty}$ such that
\begin{eqnarray}
\sup_{v \in [0,\infty[} |\lambda_n(v)-v| \rightarrow 0 \nonumber\\
\forall m \in \N\;\;\sup_{v \in [0,m]} ||x_n(\alpha_t^{-1}(\lambda_n(v)))-x(\alpha_t^{-1}(v))|| \rightarrow 0 
\label{eqSKO}
\end{eqnarray} 
Let $\gamma_n$ be defined on $[0,t[$ by $\gamma_n(u)=\alpha_t^{-1} \circ \lambda_n \circ \alpha_t(u)$. 
Then $$\sup_{u \in [0,t[}| \gamma_n(u)-u|=\sup_{u \in [0,t[} |\alpha_t^{-1}(\lambda_n(\alpha_t(u)))-u|= \sup_{v \in [0,\infty[} |\alpha_t^{-1}(\lambda_n(v))-\alpha_t^{-1}(v)|$$  As
$\alpha_t^{-1}(v)=\frac{vt}{1+v}$,  $(\alpha_t^{-1})'(v)=\frac{t}{(1+v)^2} \leq t$. It follows that 
the equations (\ref{eqSkor}) are satisfied.\\
Conversely assume that equations (\ref{eqSkor}) are satisfied. For every $m \in \N$, the derivative of $\alpha_t$ is uniformly bounded on $[0,t(1-\frac{1}{1+m})]$. It follows that $\lambda_n$ defined on $[0,\infty[$ by  $\lambda_n(v)=\alpha_t \circ \gamma_n \circ \alpha_t^{-1}(v)$ tends to $v$   uniformly on $[0,m]$ (but in general not uniformly on $[0,\infty[$). Notice also that for given $m \in \N^*$,  $\lambda_n$ is continuous strictly increasing but in general $\lambda_n([0,m]) \neq [0,m]$, i.e. $\lambda_n \notin \Lambda_m$.  In order to prove that $d_{\infty}(x_n \circ \alpha_t^{-1},x \circ \alpha_t^{-1})$ tends to $0$, we have to prove that for every $m \in \N^*$, $d_m(x_n \circ \alpha_t^{-1},x \circ \alpha_t^{-1})$ tends to $0$ (with the notations of  \cite{Bil}).   For given $m$, we have thus to replace $\lambda_n$ by some $\mu_n$ 
 in $\Lambda_m$. As $\lambda_n -Id$ goes to $0$ uniformly on $[0,m]$, for every $\epsilon>0$, one can find $n_0$ and for every $n \geq n_0$, $p_n <m$ such that $p_n \in ]m-\epsilon,m[$ and $\lambda_n(p_n) \in ]m-\epsilon,m[$. Define $\mu_n$ on $[0,m]$ by $\mu_n(v)=\lambda_n(v)$ for $v \leq p_n$, $\mu_n(m)=m$, and $\mu_n$ is affine on $[p_n,m]$. Then $\mu_n(v)$ tends to $v$ uniformly on $[0,m]$ and $\mu_n$ belongs to $\Lambda_m$.

As in  \cite{Bil} Chapter 3 Section 16, we define the function $g_m$ by:
\[g_m(t)=
\left \{\parbox{2cm}{
\begin{eqnarray} 
1\; &if& \; 0 \leq t \leq m-1 \nonumber\\
m-t &if& \; m-1 \leq t \leq m \nonumber\\\
0   &if& \; t \geq m \nonumber\
\end{eqnarray}
\nonumber
}\right.
\]
$g_m( \mu_n(v))(x_n \circ \alpha_t^{-1}(\mu_n(v)))$ tends to $g_m(v)(x \circ \alpha_t^{-1}(v))$ uniformly on $[0,p_n]$  as $n \rightarrow \infty$. 

As $p_n \in ]m-\epsilon,m[$ and $\mu_n(p_n) \in ]m-\epsilon,m[$, $|g_m(v)| \leq 
\epsilon$ and $|g_m((\mu_n(v))| \leq \epsilon$ for all $v \in ]p_n,m]$.  From equation (\ref{eqSKO}) it follows that  the sequence $x_n \circ(\alpha_t^{-1})$ is uniformly bounded on every compact subspace, thus 
$g_m( \mu_n(v))(x_n \circ \alpha_t^{-1}(\mu_n(v)))$ tends to $g_m(v)(x \circ \alpha_t^{-1}(v))$ uniformly on $[0,m]$ as $n \rightarrow \infty$.  Thus $x_n \circ \alpha_t^{-1}$ tends to $x \circ \alpha_t^{-1}$ in $D([0,\infty[,\R^d)$ for the uniform Skorokhod topology on compact spaces.  This means that $\hat d (x_n, x) \rightarrow 0$.
\hfill $\square $   
\begin{proposition}
 Let $\Omega_t=D([0,t[,\R^d)$ endowed with the uniform Skorokhod topology on compact spaces. Let $\Omega=D([0,\infty[,\R^d)$ endowed with the Skorokhod topology as defined in \cite{Bil}.
   The projection
\begin{eqnarray}
\Pi_t: D([0,\infty[,\R^d) & \rightarrow D([0,t[,\R^d) \nonumber \\
  x & \rightarrow x _{|[0,t[}
\end{eqnarray}
is  continuous.
\label{propSKO2}
\end{proposition}
{\bf Proof} Let $x_n,x \in \Omega$ such that  $x_n$ tends to $x$. 
From Theorem 16.1   of \cite{Bil},  there is a sequence $\lambda_n \in \Lambda_{\infty}$ such that
\begin{eqnarray}
\frac{\epsilon_n}{2}=\sup_{u \in [0,\infty[} |\lambda_n(u)-u| \rightarrow 0 \nonumber\\
\forall k \in \N\;\;\sup_{u \in [0,k]} ||x_n(\lambda_n(u))-x(u)|| \rightarrow 0 
\label{eqSKO2}
\end{eqnarray} 
For every $n$, $\lambda_n([0,t-\epsilon_n]) \subset [0,t[$. Let $\gamma_n$ strictly increasing from $[0,t]$ onto  $[0,t]$ such that  $(\gamma_n)_{|[0,t-\epsilon_n]}=(\lambda_n)_{|[0,t-\epsilon_n]}$,  $\gamma_n(t)=t$ and $\gamma_n$ is affine on $[t-\epsilon_n,t]$. Then 
\begin{equation}
\sup_{u \in [0,t[} |\gamma_n(u)-u| \leq \sup_{u \in [0,t-\epsilon_n]} |\lambda_n(u)-u|\leq \frac{\epsilon_n}{2} \rightarrow 0
\label{eqgamma}
\end{equation} 
For every $m>0$, there is $N $ such that for every $n \geq N$, $\epsilon_n<\frac{t}{1+m}$. Then 
\begin{eqnarray}
sup_{u \in [0,t(1-\frac{1}{1+m})]} ||x_n(\gamma_n(u))-x(u)|| & \leq & \sup_{u \in [0,t-\epsilon_n  ]} ||x_n(\lambda_n(u))-x(u)||\nonumber \\
\end{eqnarray}
From equation (\ref{eqSKO2}) it follows that 
\begin{equation}
\sup_{u \in [0,t(1-\frac{1}{1+m})]} ||x_n(\gamma_n(u))-x(u)|| \rightarrow 0
\label{eqgamma2}
\end{equation} 
From Proposition \ref{propsko}  the equations (\ref{eqgamma}) and (\ref{eqgamma2}) imply that $\hat d((x_n) _{[0,t[},x _{[0,t[}) \rightarrow 0$ in $ D([0,t[,\R^d)$.
\hfill $\square $ \\
${\cal C}_b(X)$ denotes the set of continuous bounded functions on $X$.
\begin{proposition}Let $\Omega$, and $\Omega_t$  be as in Lemma \ref{homeo}(continuous case), or as in Proposition \ref{propSKO2} (c\`adl\`ag case).  ${\cal C}_b(\Omega_t)$ is isometric to a subspace of  ${\cal C}_b(\Omega)$. Thus ${\cal C}_b(\Omega_t)$ can be identified with a subset of ${\cal C}_b(\Omega)$.
\label{propinj} 
\end{proposition} 
{\bf Proof}
Let $\Pi_t$ denote the canonical projection from $\Omega$ onto $\Omega_t$. 
From Lemma \ref{homeo} and Proposition \ref {propSKO2}, the projection $\Pi_t$ is continuous.
Thus the  map 
\begin{eqnarray}
I: {\cal C}_b(\Omega_t)& \rightarrow &{\cal C}_b(\Omega)\nonumber\\
f &\rightarrow & f \circ \Pi_t \nonumber
\end{eqnarray}
is an isometry: $||f||_{\infty}=\sup_{x \in \Omega_t}|f(x)|=\sup_{x \in \Omega}|f(\Pi_t(x)|=||f \circ \Pi_t||_{\infty}$. 
\hfill $\square $
\section{Dynamic risk measures}
\label{secdyn}                                      
In all the following, $\Omega$ is either equal to ${\cal C}_0([0,\infty[, \R^d)$ or $D_0([0,\infty[,\R^d)$. For every $t>0$, $\Omega_t$ denotes in the first case  ${\cal C}_0([0,t],\R^d)$ and in the second case  $D_0([0,t[,\R^d)$.

\begin{remark}
Denote $B_t$ the process defined on $\Omega$ by $B_t(x)=x(t)$. In case where $\Omega={\cal C}_0([0,\infty[,\R^d)$, for every $s \leq t$ $B_s$ belongs to  ${\cal C}_b(\Omega_t,\R^d)$.\\
In case where $\Omega={\cal D}_0([0,\infty[,\R^d)$, $B_t$ is continuous in $x$ if and only if $t$ is a point of continuity for $x$. Thus $B_t$ does not belong to ${\cal C}_b(\Omega)$ however $B_t$ is  measurable with respect to the $\sigma$-algebra generated by the open subsets of $\Omega_t$.
\label{remarkB}
\end{remark}
Let $c$ be a capacity   defined on ${\cal C}_b(\Omega)$.  $L^1(c)$ denotes the  Banach space obtained from completion and separation of ${\cal C}_b(\Omega)$ with respect to $c$ (cf \cite{FP} and  also \cite{BNK}). \\

\begin{proposition}
Let $c$ be a capacity on ${\cal C}_b(\Omega)$.  For every $0 \leq t \leq \infty$,  ${\cal C}_b(\Omega_t) \subset L^1(c)$ (with the notation $\Omega_{\infty}=\Omega$). Denote $L^1_t(c)$ the closure of  ${\cal C}_b(\Omega_t)$ in the Banach space $L^1(c)$.
For every $s \leq t$, $L^1_s(c) \subset L^1_t(c)$.
\label{propcyl}
\end{proposition}
{\bf Proof} For every $t>0$, it follows from Proposition \ref{propinj} that ${\cal C}_b(\Omega_t) \subset {\cal C}_b(\Omega)$. As in Proposition \ref{propSKO2}, one can prove for $s \leq t$ that the map $x \rightarrow x_{|[0,s[}$ induces a continuous projection from $D([0,t[,\R^d)$ onto $D([0,s[,\R^d)$. Thus 
${\cal C}_b(\Omega_s)\subset {\cal C}_b(\Omega_t)$. 
 and $L^1_s(c) \subset L^1_t(c)$.
\hfill $\square $.\\
In all the following   $\Omega_{\infty}$ means $\Omega$ and $L^1_{\infty}(c)$ means $L^1(c)$. 
\begin{definition}
A dynamic risk measure on $L^1(c)$ is a family $(\rho_{st})_{\{0 \leq s \leq t \leq \infty\}}$ 
where $\rho_{st}: L^1_t(c)\rightarrow L^1_s(c)$ satisfies 
\begin{enumerate}
\item monotonicity: $\rho_{st}(X) \geq \rho_{st}(Y)$ for all $X,Y \in L^1_t(c)$ such that $X \leq Y$.
The above inequalities are  with respect to the order defined on $L^1(c)$ (cf Section 2.1 of \cite{BNK}).
\item translation invariance: $\forall X \in  L^1_t(c),\; \forall Y \in  L^1_s(c)$, 
$$\rho_{st}(X+Y)=\rho_{st}(X)-Y$$
\item convexity:  $\forall X,Y  \in  L^1_t(c)$, $\forall \lambda \in [0,1]$, 
 $$\rho_{st}(\lambda X+ (1-\lambda) Y)\leq \lambda \rho_{st}(X)+ (1-\lambda)\rho_{st}(Y)$$
\end{enumerate}
\label{defdyn}
\end{definition}
Before presenting our results on  dynamic risk measures, we recall some results obtained in a previous paper \cite{BNK}.
 In this paper, we have proved that every continuous linear form on $L^1(c)$  is represented by a bounded signed measure on $(\Omega,{\cal B}(\Omega))$ and that $K_+$,  the non negative part of the unit ball of the dual of $L^1(c)$ is metrizable compact for the weak* topology.
Assume now that $c(f)=\sup_{P \in {\cal P}}E_p(|f|^p)^{\frac{1}{p}}$ where ${\cal P}$ is a weakly relatively compact set of probability measures. We have proved  that even if this set is non dominated, there is a  countable subset  $\{ Q_n, n \in \N \}$ of $\mathcal{P}$ such that, for every $X\in L^1(c)$,
\begin{equation}
c(X)=\sup_{n\in\N} \left( E_{Q_n}\left( |X|^p\right)\right)^{\frac{1}{p}}
\label{eqc}
\end{equation}
(cf Theorem 4.1 of \cite{BNK}
 for more details and the proof).\\
 A very important result in \cite{BNK} is the existence of a canonical equivalence class of measures $\mu$ in $K_+$ such that an element  of $L^1(c)$ is completely determined if one  knows it $\mu$ a. s. This equivalence class is referred to as the canonical $c$-class. It is characterized by the following property:  $\mu \in K_+$ belongs to the canonical $c$-class if and only if  
$$\forall X, Y \in L^1(c), \{X =Y \;\;\;\mu \; a.s.\} \;\Longleftrightarrow\;\{ X=Y\; \text{in}\; L^1(c)\}$$ 
In all the following $P$ is a given probability measure belonging to the canonical $c$-class. A particular choice is $  \sum_{n \in \N} \frac{1}{2^{n+1}}Q_n$ where $\{Q_n, \; n \in \N\}$ satisfies (\ref{eqc}).
We  have also proved (see Theorem 3.1 of  \cite{BNK}
) that every convex risk measure $\rho$ on $L_1(c)$ is continuous and admits the following representation: 
\begin{equation}
\forall X\in L^1(c),\ \rho\left( X\right)=\sup_{Q\in\mathcal{P'}}(
E_{Q}[-X]-\alpha\left( Q\right))
\label{eqmag}
\end{equation}
where \begin{equation}
\alpha\left(Q\right)=\sup_{X\in L^{1}(c)}\left( E_{Q}[-X]-\rho\left( X\right)\right)
\label{eqma}
\end{equation}
and $\mathcal{P'}$ is the set of probability measures on $(\Omega,{\cal B}(\Omega))$ belonging to $L^1(c)^*$.\\
\label{thmrep1}
 \section{Representation of Dynamic Risk Measures}
\label{secrepdyn}
In this section we  want to prove a dual representation result  for dynamic risk measures on $L^1(c)$.
${\cal P}$ is  a weakly relatively compact set of probability measures on $\Omega$ and $c$ is the capacity $c(f)=\sup_{P \in {\cal P}}E_P(|f|^p)^{\frac{1}{p}}$. In all this section, $({\cal B}_t)_{t \in \R^+}$ denotes  the right continuous filtration  generated by the open sets of $\Omega_t$  ($({\cal B}_t)_{t \in \R^+}$ is not completed). $P$ is a given probability measure in $K_+$ (the unit ball of the dual of $L^1(c)$) belonging to the canonical $c$-class.
\begin{lemma} 
 $L^1_t(c) \subset L^1(\Omega,{\cal B}_t,  P)$. Furthermore every element of $L^1_t(c)$ is characterized by its class in $L^1(\Omega,{\cal B}_t,P)$.
\label{lemmadomine}
\end{lemma} 
 {\bf Proof}
Every element of $L^1_t(c)$ is limit for the $c$-norm of a sequence of continuous functions on $\Omega_t$. It  follows that every element of $L^1_t(c)$ can be represented by a $(\Omega,{\cal B}_t)$ measurable function $f$. 
 As $P$ belongs to  $ K_+$, for  every $f \in L^1(c)$, $E_{ P}(|f|) \leq c(f)$.\\
Furthermore, let $f,g \in L^1_t(c)$ such that $f=g\;P\;a.s.$. Then $|f-g|=0\;P\;a.s.$. From Theorem 4.2 of \cite{BNK}, it follows that $f=g$ in $L^1_t(c)$.
\hfill $\square $ \\
In particular for every $f \in L^1_s(c)$, $E_{ P}(|f|)< \infty$. And this allows to compose  $E_{ P}$ with $\rho_{st}$. Denote $\tilde \rho_{st}=E_{P}o\rho_{st}$

\begin{definition}
Denote ${\cal M}^+_{s,t}(P)$ the set of non negative continuous linear forms on $L^1_t(c)$ represented by a probability measure on $(\Omega, {\cal B}(\Omega_t))$ whose restriction to ${\cal B}_s$ is absolutely continuous with respect to $P$.\\
Define 
\begin{itemize}
\item the acceptance set 
$${\cal A}_{st}=\{X \in L^1_{t}(c) \; |\; \rho_{st}(X) \leq 0 \; \;\}$$ The above inequality is with respect to the order defined on $L^1_s(c)$ (cf Section 2.1 of \cite{BNK}. As $P$ belongs to the canonical $c$-class, $ \rho_{st}(X) \leq 0$ is equivalent to $ \rho_{st}(X) \leq 0\; P \; a.s.$
\item Denote
\begin{equation}
{\tilde{\cal A}}_{st}=\{\sum_{i \in I} X_i 1_{A_i}\;|X_i \in {\cal A}_{st},\;\; A_i \in {\cal B}_s\;  (A_i)_{i \in I}\; finite \;  partition \;of\; \Omega\}
\label{Atilde}
\end{equation}
\item  the minimal penalty for every probability measure $Q \in {\cal M}^+_{s,t}(P)$
$$\alpha^m_{st}(Q)=Q\esssup_{X \in {\cal A}_{st}}E_Q(-X |{\cal B}_s)\;\; Q \; a.s.$$
\item the $Q$-acceptance set for every probability measure  $Q \in {\cal M}^+_{s,t}(P)$
\begin{equation}
{\cal A}_{st}(Q)=\{X \in L^1_{t}(c) \; |\; \rho_{st}(X) \leq 0 \; \;Q\; a.s.\}
\label{eqAQ}
\end{equation}
\end{itemize}
\end{definition}

\begin{lemma}
Let $Q$ be a probability measure probability measure in ${\cal M}^+_{s,t}(P)$.
\begin{enumerate}
\item The minimal penalty admits the following representations:
\begin{eqnarray}
\label{eq_6} 
 \alpha^m_{s,t}(Q) & = & Q \esssup_{X \in {\cal A}_{s,t}} E_Q(-X| {\cal B}_s) \nonumber \\ 
 & = & Q \esssup_{X \in {\cal A}_{s,t}(Q)} E_Q(-X| {\cal B}_s) \nonumber \\ 
 & = &Q \esssup_{X \in L^1_t(c)} 
(E_Q(-X| {\cal B}_s)-\rho_{s,t}(X)) \label{eqpen}\\ 
& = &Q  \esssup_{X \in \tilde{\cal A}_{s,t}}E_Q(-X| {\cal B}_s). \nonumber 
\end{eqnarray}
\item
${\cal Z}=\{ E_Q(-X| {\cal B}_s);\;X \in \tilde {\cal A}_{s,t}\}$
is a lattice upward directed in $L^1(\Omega,{\cal B}_s,Q)$.
\end{enumerate}
\label{lemmalat}
\end{lemma}
{\bf Proof}
\begin{enumerate}
\item
Let $Y$ in $L_s^1(c)$. From Lemma 4.3 of \cite{BNK},  $Y \leq 0$ for the order in $L_s^1(c)$ iff $Y \leq 0\;P \;a.s.$. It follows that for every $Q \in {\cal M}^+_{s,t}(P)$, ${\cal A}_{s,t} \subset  {\cal A}_{s,t}(Q) \subset L^1_t(c)$ and 
\begin{eqnarray} Q \esssup_{X \in {\cal A}_{s,t}} E_Q(-X| {\cal B}_s) & \leq & Q \esssup_{X \in {\cal A}_{s,t}(Q)} E_Q(-X| {\cal B}_s)\nonumber \\
 & \leq & Q \esssup_{X \in L^1_t (c)} (E_Q(-X| {\cal B}_s)-\rho_{s,t}(X)) \nonumber
\end{eqnarray}
On the other hand, for every $X \in  L^1_t(c)$, $X +\rho_{st}(X) \in {\cal A}_{s,t}$, so the above inequalities are in fact equalities.\\
${\cal A}_{s,t} \subset  \tilde {\cal A}_{s,t}$ thus 
$Q \esssup_{X \in {\cal A}_{s,t}} E_Q(-X| {\cal B}_s) \leq Q \esssup_{X \in \tilde{\cal A}_{s,t}}E_Q(-X| {\cal B}_s)$. Let $f$ be  ${\cal B}_s$-measurable such that for all X in ${\cal A}_{s,t}$, $f \geq  E_Q(-X| {\cal B}_s) \; Q \; a.s.$. For every $A \in {\cal B}_s$, $f1_A \geq  E_Q(-X| {\cal B}_s)1_A= E_Q(-X1_A| {\cal B}_s)$. It follows that for every $ (A_i)_{i \in I} \in {\cal B}_s$ finite partition of $ \Omega$, for every family $X _i \in{\cal A}_{s,t}$, 
$$f = \sum_{i \in I}f 1_{A_i}\geq  E_Q(-\sum _{i \in I}X_i 1_{A_i}|{\cal B}_s)$$
And thus $f  \geq Q \esssup_{X \in \tilde{\cal A}_{s,t}}E_Q(-X| {\cal B}_s)$.
This proves the equations (\ref{eqpen}).
\item Let $X, Y \in \tilde{\cal A}_{s,t}$. Let $A= \{\omega\; |\;E_Q(-X| {\cal B}_s) \geq E_Q(-Y| {\cal B}_s)\}$. $A$ is ${\cal B}_s$ measurable. $Z=X1_A+Y1_{\Omega -A} \in \tilde{\cal A}_{s,t}$, and $E_Q(-Z |{\cal B}_s)= \sup (E_Q(-X| {\cal B}_s), E_Q(-Y| {\cal B}_s))$.\hfill $\square $
\end{enumerate}

In all the following for $X$ in $L^1(c)$, the equality $X=Y \:P \;a.s.$  means that $X$ is the unique element of $L^1(c)$ equal to $Y \;\; P$ a.s. \\
 For $0 \leq s \leq t \leq \infty$, denote ${\cal Q}_{s,t}(P)$  the set of probability measures in the dual of $L^1_t(c)$  whose restriction to ${\cal B}_s$ is equal to $P$.
\begin{theorem}
 Let $\rho_{st}$ be a dynamic risk measure on $L^1(c)$.   The following conditions are equivalent:
\begin{enumerate}
\item The dynamic risk measure is strongly convex i.e. satisfies:
$\forall X, Y \in L^1_t(c)\;\;\forall f \in {\cal C}_b(\Omega_{s}),\;\; 0 \leq f \leq 1,$ 
$$\rho_{st}(Xf+Y(1-f)) \leq \rho_{st}(X)f+\rho_{st}(Y)(1-f)$$
\item The acceptance set ${\cal A}_{st}$ is strongly convex i.e.
$\forall X, Y \in {\cal A}_{st}\;\;\forall f \in {\cal C}_b(\Omega_{s}), 0 \leq f \leq 1,$ $$Xf+Y(1-f) \in {\cal A}_{st}$$
\item There is a countable set $\{Q_n, n \in \N\}$ of probability measures in ${\cal Q}_{s,t}(P)$ such that 
$\rho_{st}$ admits the following dual representation.
\begin{equation}
X\in L_t^1(c),\ \rho_{st}(X)=P \esssup_{ n \in \N} (E_{Q_n}(-X|{\cal B}_s)- \alpha^m_{st}(Q_n))\;\;P \;a.s.
\label{eqrep20}
\end{equation}
\item
$\rho_{st}$ admits the following dual representation.
\begin{equation}
\forall X\in L_t^1(c),\ \rho_{st}(X)=P \esssup_{Q \in {\cal Q}_{s,t}(P)} (E_Q(-X|{\cal B}_s)- \alpha^m_{st}(Q))\;\;P \;a.s.
\label{eqrep21}
\end{equation}
\end{enumerate} 
 \label{thm3}
\end{theorem}
From Section 4.2. of \cite{BNK}, an element of $L^1(c)$ is fully determined by its expresssion $P$ a.s.. Thus the representation (\ref{eqrep20}) or (\ref{eqrep21}) characterizes completely the element $\rho_{st}(X)$ of $L^1_t(c)$.\\ Notice that the set $\tilde{\cal A}_{s,t}$ has been introduced in order to have the  lattice property expounded in Lemma \ref{lemmalat}. Indeed in the general case  $L^1_t(c)$ does not contain the characteristic functions of ${\cal B}_s$ measurable sets. Thus  the lattice property is not satisfied for $\{ E_Q(-X| {\cal B}_s);\;X \in  {\cal A}_{s,t}\}$.\\
{\bf Proof}
This proof is inspired by the proof of Proposition 1 of \cite{JBN}. However it is more complicated here due to the fact that $L^1(c)$ does not contain in general all the bounded ${\cal B}(\Omega)$-measurable functions.\\
The implications ${\it 4}. \Longrightarrow {\it 1}.$ and ${\it 1}. \Longrightarrow {\it 2}.$ are easily verified. ${\it 3}. \Longrightarrow {\it 4}.$ follows easily from equation (\ref{eqpen}).\\
 Now prove that 
 ${\it 2}.$ implies ${\it 3}$.
In order to prove the representation (\ref{eqrep20}) it is not restrictive to assume that $\rho_{st}$ is normalized, i.e. satisfies $\rho_{st}(0)=0$. From the theorem of representation of convex risk measures on $L^1(c)$ (Theorem 3.3 of \cite{BNK}),  there is a countable set ${\cal Q}_0=\{Q_n,\;n \in \N\}$ of probability measures on $(\Omega,{\cal B}(\Omega))$ belonging to the dual of $L_t^1(c)$ such that   $\tilde \rho_{st}=E_{ P}o \rho_{st}$ admits the representation
\begin{equation}
\tilde \rho_{st}(X)= \sup_{Q \in {\cal Q}_0} (E_Q(-X)-\alpha (Q))
\label{eqrep3}
\end{equation}
where $\alpha(Q)=\sup_{X \in L^1_t(c)} (E_Q(-X)-\tilde \rho_{st}(X))$. 
Of course one can assume that for every $Q \in {\cal Q}_0$, $\alpha(Q)$ is finite. 
 We adapt the proof of Proposition 1 of \cite{JBN}.\\
{\it i)} First step: We prove that for every $Q \in {\cal Q}_0$ such that $\alpha (Q)$ is finite,  the restriction of $Q$ to ${\cal B}_s$ is equal to $P$.\\
 In the proof of Proposition 1 of \cite{JBN}, this property was obtained computing $\rho_{st}(\beta 1_A)$  for $\beta \in \R$ and $A \in {\cal B}_s$. This cannot be done here because in general, $1_A \notin L^1_t(c)$. \\
We have to use here the regularity of the measures $P$ and $Q$ which follow from Theorem 1.1 of \cite{Bil}. Thus for every ${\cal B}_s$ measurable set $A$ there is a family of bounded continuous functions $f_n$ on $\Omega_s$, $0 \leq f_n \leq 1$ such that $f_n$ tends to $1_A$ $ P+Q$ a.s. From normalization and translation invariance of $\rho_{st}$, the restriction of $\rho_{st}$ to $L^1_s(c)$ is equal to $-Id_{L^1_s(c)}$. Thus for every $n$ and every $\beta \in \R$, $\rho_{st}(\beta f_n) = -\beta f_n$. From the representation of $\tilde \rho_{st}$ (equation (\ref{eqrep3})), it follows that 
\begin{equation}
\alpha(Q) \geq - E_Q(\beta f_n)+E_{P}(\beta f_n)
\label{eqeg}
\end{equation}
As $1_A$ is the limit of $f_n$ $P+Q$ a.s., using the dominated convergence Theorem, one can pass to the limit in (\ref{eqeg}) and get that    
$\alpha(Q) \geq \beta (P(A)-Q(A))$ for every $\beta \in \R$. As $\alpha(Q) <\infty$, necessarily, $Q(A)=P(A)$ for every $A \in {\cal B}_s$. We have thus proved that ${\cal Q}_0 \subset {\cal Q}_{s,t}(P)$.\\
{\it ii)} The next step  consists in proving that for every $Q$ such that $\alpha(Q) < \infty$, $\alpha(Q)= E_Q(\alpha^m_{st}(Q))$.\\ By definition 
$\alpha(Q)=\sup_{X \in L^1_t(c)}(E_Q(-X)-\tilde \rho_{st}(X))$. From i), it follows that $\alpha(Q)=\sup_{X \in L^1_t(c)}(E_Q(-X- \rho_{st}(X)))$. The expression of $\alpha^m_{st}(Q)$ is given by (\ref{eqpen}), thus
\begin{equation} 
\alpha(Q)\leq  E_Q(\alpha^m_{st}(Q))
\label{eqmin}
\end{equation}
 From  the expression of $\alpha^m_{st}(Q)$, the lattice property of ${\cal Z}$( cf  Lemma \ref{lemmalat}) , and Lemma 1 of \cite{JBN}, it follows that $E_Q(\alpha^m_{st}(Q))=\sup_{X \in \tilde{\cal A}_{s,t}}E_Q(-X)$. Let $\epsilon>0$, there is an element $Z \in \tilde {\cal A}_{st}$
 such that  $E_Q(-Z)> E_Q(\alpha^m_{st}(Q))- \epsilon$. $Z\in \tilde {\cal A}_{st}$ can be written as a finite sum
$Z=\sum_{i \in I} X_{i} 1_{A_{i}}$ $X_{i} \in {\cal A}_{st},\;\; (A_{i})_{i \in I} \in {\cal B}_s$ is a finite  partition of $\Omega$. It follows from the regularity of $P$ that  there is a  sequence $(f^k_{i})_{k \in \N, i \in I}$ of continuous  functions  $0 \leq f^k_{i} \leq 1$, $ \sum_{i \in I}  f^k_{i}=1$, such that for every $i \in I$, $1_{A_{i}}= \lim_{k \rightarrow \infty} f^k_{i}\;\;  P \;a.s.$ and  thus also $ Q \;a.s.$. Furthermore every $Q$ in the representation of $\tilde  \rho_{st}$ is a continuous linear form for the $c$ norm. Thus there is a constant $K$ (depending on $Q$) such that for every $X \in L^1_t(c)$, $|E_Q(X)| \leq Kc(X)$. As $c(X)=c(|X|)$, we get $E_Q(|X|) \leq K c(X)$. As $0 \leq f^k_i \leq 1$, $|\sum_{i \in I} X_{i} f^k_{i}| \leq  \sum_{i \in I} |X_{i}|$ and $E_Q(\sum_{i \in I} |X_{i}|) \leq K \sum_{i \in I} c(|X_{i}|) < \infty$. Thus we can apply the dominated convergence theorem,  it follows that there is $k$ such that  $E_Q(-\sum_{i \in I} X_{i} f^k_{i})> E_Q(\alpha^m_{st}(Q))- 2\epsilon$. By hypothesis 2., $X=\sum_{i \in I} X_{i} f^k_{i} \in {\cal A}_{st}$. As $\alpha(Q) \geq \sup_{X \in {\cal A}_{st}} E_Q(-X)$ it follows that  $\alpha(Q) \geq E_Q(\alpha^m_{st}(Q))$. Thus from (\ref{eqmin}), for every  $ Q$  such that $\alpha(Q)<\infty$,
\begin{equation}
 \alpha(Q) = E_Q(\alpha^m_{st}(Q))
\label{eqpen2}
\end{equation}
{\it iii)} last step\\
The inequality 
\begin{equation}
\rho_{st}(X) \geq \esssup_{Q \in {\cal Q}_0} E_Q(-X|{\cal B}_s)- \alpha^m_{st}(Q)\;\; P\;a.s.
\label{eqrep4}
\end{equation}
follows from equality (\ref{eqpen}).
From the equations (\ref{eqrep3}) and (\ref{eqpen2}), it follows that 
$$\tilde \rho_{st}(X)= \sup_{Q \in {\cal Q}_0} E_Q(-X-\alpha^m_{st}(Q))$$
Furthermore the restriction of every $Q \in {\cal Q}_0$ to ${\cal B}_s$ is equal to $P$. So
$$\tilde \rho_{st}(X)= \sup_{Q \in {\cal Q}_0} E_{P}(E_Q(-X|{\cal B}_s)-\alpha^m_{st} (Q))$$
It follows that 
\begin{equation} 
E_{P}(\rho_{st}(X))\leq E_{P} (\esssup_{Q \in {\cal Q}_0}(E_Q(-X|{\cal B}_s)-\alpha^m_{st} (Q))
\label{eqrep6}
\end{equation}
From the inequalities (\ref{eqrep4}) and (\ref{eqrep6}), as ${\cal Q}_0$ is a countable subset of ${\cal Q}_{st}(P)$,  we get the announced representation (\ref{eqrep20}).
\hfill $\square $ 
\begin{corollary} Assume that the risk measure $\tilde \rho_{s,t}$ is normalized and majorized by a  sublinear risk measure. Then for every $X \in L^1_t(c)$, there is a probability measure $Q_X$  such that 
\begin{eqnarray}
\rho_{st}(X) & = &  E_{Q_X}(-X|{\cal B}_s)- \alpha^m_{st}(Q_X)\; \; P \;a.s.\nonumber\\
 &  = & \esssup_{Q \in {\cal Q}_{s,t}(P)} (E_Q(-X|{\cal B}_s)- \alpha^m_{st}(Q))\; \; P \;a.s. \nonumber\\
\end{eqnarray}
\label{c1}
\end{corollary}
{\bf Proof}
From Proposition 3.1 and Theorem 3.2 of \cite{BNK}, there is a probability measure $Q_X$  such that 
\begin{equation}
\tilde \rho_{st}(X)= E_{Q_X}(-X)- \alpha(Q_X)
\label{eq23}
\end{equation}
 $\alpha(Q_X) < \infty$ so from the proof of the preceding theorem, $Q_X$ is in  ${\cal Q}_{s,t}(P)$. Thus from (\ref{eqpen2}) and (\ref{eq23}) and (\ref{eqrep4}), it follows that $\rho_{st}(X)  =  E_{Q_X}(-X|{\cal B}_s)- \alpha^m_{st}(Q_X)$
\hfill $\square $ 

\section{Continuity of Dynamic risk measures} 

 In order to prove continuity results for dynamic risk measures, we give first an expression of the $c$ norm for non negative elements of $L^1(c)$ as a supremum of linear forms. Denote $q$ the conjugate exponent of $p$ ($\frac{1}{p}+ \frac{1}{q}=1$). From Theorem 4.1 of \cite{BNK}
 recalled in Section \ref{secdyn}, there exists a countable set $(Q_n)_{n \in \N}$   of probability measures such that 
$$\forall X \in { L}^p(c),\; c(X)=\sup_{n \in \N} (E_{Q_n}(|X|^p))^{\frac{1}{p}}$$
\begin{lemma} Denote $L=\{g, \; {\cal B}_s\; measurable\;|\sup_{n \in \N}E_{Q_n}(|g|^q)^{\frac{1}{q}} \leq 1\}$, and $L_+=\{g \in L|g \geq 0\}$. For every $s>0$, for every $X  $ in $L^1_s(c)$, 
\begin{equation}
c(X)=\sup_{n \in \N}(E_{Q_n}(|X|^p)^{\frac{1}{p}}=\sup_{n\in \N}(\sup_{\{g \in L_+\}}(E_{Q_n}(|X|g)))
\label{eqN1}
\end{equation}
Furthermore there is $g_0 \in L_+$ such that $c(X)=\sup_{n\in \N}E_{Q_n}(|X|g_0)$
\label{norm1}
\end{lemma}

{\bf Proof}
We already know that  $c(X)=\sup_{n \in \N}(E_{Q_n}(|X|^p)^{\frac{1}{p}}$ for every $X \in L^1_s(c)$. The inequality 
\begin{equation}
c(X) \geq \sup_{n\in \N}(\sup_{\{g \in L_+\}}(E_{Q_n}(|X|g)))
\label{eqN2} 
\end{equation}
follows then from H\"older inequality. Let  $X \in L^1_s(c)$, $c(X)<\infty$. Denote $g_0=\frac{|X|^{\frac{p}{q}}}{c(X)^{p-1}}$, $g_0 \geq 0$.
For every $Q_n$, $E_{Q_n}(g_0^q) = E_{Q_n}(\frac{|X|^p}{c(X)^p})\leq 1$, thus $g_0 \in L_+$. For every $n$, 
$E_{Q_n}(|X|g_0)=E_{Q_n}(\frac{|X|^p}{c(X)^{p-1}})$, thus $\sup_{n \in \N}E_{Q_n}(|X|g_0)=c(X)$. This proves the result.
\hfill $\square $ \\
Define $\phi$  on $L^1_s(c)$ by 
\begin{equation}
\phi(X)=\sup_{n \in N}\sup_{\{Y \in L_+\}}E_{Q_n}(XY)
\label{eqro}
\end{equation}
\begin{theorem} For every $s \leq t$, the conditional risk measure $\rho_{st}$ defined on $L^{1}_t(c)$ with values in $L^1_s(c)$ is continuous for the $c$ norm.
\label{thm6}
\end{theorem}
{\bf Proof}
We want to prove that every $\rho_{st}$ satisfying the axioms of Definition \ref{defdyn} is continuous for the $c$ norm. For every $X \in L^1_t(c)$ consider the function $\rho_{X,s,t}$ defined on $L^1_t(c)$ by $\rho_{X,s,t}(Y)=\rho_{s,t}(X+Y)$. $\rho_{X,s,t}$ satisfies also the axioms of Definition \ref{defdyn} and the continuity of $\rho_{X,s,t}$ in $0$ is equivalent to the continuity of $\rho_{s,t}$ in $X$. Thus it is enough to prove that every conditional risk measure $\rho_{s,t}$ is continuous in $0$. Notice also that considering $\bar \rho_{st}$ defined by $\bar \rho_{st}(X)=\rho_{s,t}(X)-\rho_{s,t}(0)$ we can furthermore assume that $\rho_{st}(0)=0$\\
Using the definition of $\phi$, one can remark that $\phi$ is convex, and for $X \leq Y$, $\phi(X) \leq \phi(Y)$.
Let $\tilde \phi$ be the function defined on $L^1_t(c)$ by $\tilde \phi(X)=\phi(\rho_{st}(X))$. Then, $\tilde \phi$ obviously satisfies the axioms of convexity and monotonicity. Thus, thanks to Theorem 1 of \cite{BF}, it is continuous for the $c$ norm. Notice that  $\tilde \phi$ is not necessarily a risk measure as it does not necessarily satisfy the translation invariance. $\tilde\phi$ is continuous in $0$ thus 
for every 
$ \epsilon>0$, there is $\eta>0$  such that $c(X)<\eta$ implies 
$\tilde\phi(X)< \epsilon$. For every $X\leq 0$, $X \in L^1_t(c)$, $\rho_{st}(X) \geq 0$ and from Lemma \ref{norm1}, $c(\rho_{st}(X))= \phi(\rho_{st}(X))=\tilde \phi(X)$. Thus for every $X \in L^1_t(c)$,
 $c(X)<\eta$ implies 
$c(\rho_{st}(-|X|))< \epsilon$.\\
As $\rho_{st}$ is convex and satisfies $\rho_{st}(0)=0$, it follows that $|\rho_{st}(X)| \leq \rho_{st}(-|X|)$. It follows that for all $X \in L^1(c)$ such that $c(X) < \eta$, $c(\rho_{st}(X))< \epsilon$.
This proves the continuity of $\rho_{st}$ for the $c$ norm at point zero and thus everywhere.
\hfill $\square $

\section{Time consistency and regularity of paths}
The goal of this section is to prove in the context of model uncertainty two results which were proved for dynamic risk measures on $L^{\infty}$ spaces. First as in \cite{BN02} we  want to characterize the time consistency property by a cocycle condition on the penalty. Second we want to prove the regularity of paths for convex time consistent dynamic risk measures on $L^1(c)$. This last result was first proved for sublinear i.e. coherent time consistent dynamic risk measures on $L^{\infty}$ spaces by Delbaen \cite{D}. It was extended to convex time consistent dynamic risk measures on $L^{\infty}$ spaces in \cite{JBN}. Notice however that the regularity of paths in the present case cannot be deduced from the result on $L^{\infty}$ spaces. Indeed as already mentioned $L^1(c)$ does not contain all the bounded Borelian functions,  furthermore a dynamic risk measure on $L^1(c)$ cannot be always extended to $L^{\infty}(P)$ as the probability measures in the dual representation of $\rho_{st}$  are not all absolutely continuous with respect to $P$ (Theorem \ref{thm3}). Thus we have to give new proofs. \\
 The notations are those of Section \ref{secrepdyn}. ${\cal B}_t$ is the right continuous filtration generated by the Borelian sets of $\Omega_t$.   $P$ is  a given probability measure in $K_+$ belonging to the canonical $c$-class.
${\cal M}^+_{s,t}(P)$ (resp. 
${\cal Q}_{s,t}(P)$) denotes the set of non negative continuous linear forms on $L^1(c)$ represented by a probability measure on $(\Omega, {\cal B}_t)$ whose restriction to ${\cal B}_s$ is absolutely continuous with respect to $P$ (resp.
 equal to $P$).
For $Q$ in ${\cal M}^+_{s,t}(P)$, the penalty $\alpha_{s,t}^m\left( Q\right)$ is given by the formula (\ref{eq_6}), and the $Q$-acceptance set is given by equation (\ref{eqAQ}).

\subsection{Characterization of time consistency}

\begin{theorem}
Assume that $\left( \rho_{s,t}\right)_{s,t}$  is a dynamic risk measure on $L^1(c)$ which satisfies the equivalent conditions of Theorem \ref{thm3}. Let $r\leqslant s\leqslant t$ be three instants of  time. The following conditions are equivalent:
\begin{itemize}
\item[i)] $\forall X\in L_{t}^p,\ \rho_{r,t}\left( X\right)=\rho_{r,s}\left( -\rho_{s,t}\left( X\right)\right)$
\item[ii)] For every probability measure $Q$ in ${\cal M}^+_{s,t}(P)$,  
$${\mathcal{A}}_{r,t}\left( Q\right)={\mathcal{A}}_{r,s}\left( Q\right)+{\mathcal{A}}_{s,t}$$
\item[iii)] For every probability measure $Q$ in ${\cal M}^+_{s,t}(P)$
$$\alpha_{r,t}^m\left( Q\right)=\alpha_{r,s}^m\left( Q\right)+E_Q\left( \alpha_{s,t}^m\left( Q\right)\vert \mathcal{B}_r\right)\ Q \;a.s.$$
\end{itemize}
\label{thmtc}
\end{theorem}
{\bf Proof}\ \\
The proof of $i) \Rightarrow ii)$ is the same as the proof of Theorem 1 in \cite{JBN}.\\
 $iii) \Rightarrow i)$. From the theorem of representation, Theorem \ref{thm3}, equation(\ref{eqrep20}), there is a countable set $\{Q_n,\; n \in \N\}$ of probability measures in  ${\cal Q}_{s,t}(P)$ such that 
$$\rho_{st}(X)=P \esssup_{ n \in \N} (E_{Q_n}(-X|{\cal B}_s)- \alpha^m_{st}(Q_n))$$
Notice that $Q_n$  is a probability measure on $\left( \Omega, \mathcal{B}_t\right)$ whose  restriction  to $\mathcal{B}_s$ is equal to $P$. $Q_n$   is not in general absolutely continuous with respect to $P$. Therefore, we introduce the probability measure: $\tilde P= \frac{P}{2} + \sum_{n \in \N}\frac{Q_n}{2^{n+2}}$. Denote ${\cal Q}$ the set of probability measures in the dual of $L^1_t(c)$ absolutely continuous with respect to $\tilde P$ whose restriction to ${\cal B}_s$ is equal to $\tilde P$ (i.e. $P$). The set 
$\{E_{R}(-X|{\cal B}_s)- \alpha^m_{st}(R),\; R \in {\cal Q}\}$ is then a lattice upward directed.\\
 We  can then apply similar  arguments to those used  in the proof of Theorem 1 in \cite{JBN} to prove that $iii) \Rightarrow i)$.
\\
Let's now prove $ii) \Rightarrow iii)$. From equations (\ref{eq_6}), and hypothesis $ii)$,
\begin{eqnarray*}
\alpha_{r,t}^m\left( Q\right)&=&\makebox{ess sup}_{Y\in\mathcal{A}_{r,s}(Q)}E_Q\left[ -Y\vert \mathcal{B}_r\right]+\makebox{ess sup}_{Z\in\mathcal{A}_{s,t}}E_Q\left[ -Z\vert \mathcal{B}_r\right]\;Q \; a.s.\\
&=&\alpha_{r,s}^m\left( Q\right)+\makebox{ess sup}_{Z\in\mathcal{A}_{s,t}}E_Q\left[ -Z\vert \mathcal{B}_r\right] \;Q \; a.s.
\end{eqnarray*}
It only remains to prove that 
\begin{equation}
\makebox{ess sup}_{Z\in\mathcal{A}_{s,t}}E_Q\left[ -Z\vert \mathcal{B}_r\right]=E_Q\left[ \alpha_{s,t}^m\left( Q\right) \vert  \mathcal{B}_r\right]\label{eqE}
\end{equation} 
Using the formula (\ref{eqpen}), the first inequality,
$$\makebox{ess sup}_{Z\in\mathcal{A}_{s,t}}E_Q\left[ -Z\vert \mathcal{B}_r\right]\leqslant E_Q\left[ \alpha_{s,t}^m\left( Q\right) \vert  \mathcal{B}_r\right]$$
is obvious. \\
In order to prove the converse inequality, we recall (see Lemma 3) that $\left\{ E_Q\left[-X\vert \mathcal{B}_s \right]; X\in\widetilde{\mathcal{A}}_{s,t}\right\}$ is a lattice upward directed. Therefore, using formula $(\ref{eqpen})$, we deduce from Proposition VI-1-1 of  \cite{nev}  the existence of a  sequence $\left( Z_n\right)_{n\in\N}$ in $\widetilde{\mathcal{A}}_{s,t}$ such that
$\alpha_{s,t}^m\left( Q\right)$ is $Q$-a.s. the increasing limit of $\left( E_Q\left[-Z_n\vert \mathcal{B}_s \right] \right)_{n\in\N}$.\\
$E_Q\left[\alpha_{s,t}^m\left( Q\right)\right]$ is  the increasing limit of $\left( E_Q\left[-Z_n \right] \right)_{n\in\N}$. 
\begin{itemize}
\item
If $E_Q\left[\alpha_{s,t}^m\left( Q\right)\right]<+\infty$, let $\epsilon>0$, there exists $n\in\N$ such that
$$ E_Q\left[-Z_n \right] \leq E_Q\left[\alpha_{s,t}^m\left( Q\right)\right]\leqslant E_Q\left[-Z_n \right] +\epsilon$$
Now, as $\left( Z_n\right)_{n\in\N}\in \widetilde{\mathcal{A}}_{s,t}$, $Z_n=\sum_{i\in I}X_{i,n} \mathbf{1}_{A_{i,n}}$ where $X_{i,n} \in {\mathcal{A}}_{s,t}$, $A_{i,n}\in \mathcal{B}_s$ and $\left( A_{i,n}\right)_{i\in I}$ is a finite partition of $\Omega$ in  $\mathcal{B}_s$.
As in the proof of Theorem \ref{thm3}, there exists a sequence of bounded continuous functions $f_{i,n}^k$ on $\Omega_{[0,s]}$, $0\leqslant f_{i,n}^k \leqslant 1$ such that $f_{i,n}^k$ tends to $\mathbf{1}_{A_{i,n}}$ $Q$ a.s. when $k$ tends to $+\infty$.\\
Let $n$ fixed. Notice that  $\left|\sum_{i\in I}X_{i,n} f_{i,n}^k\right| \leqslant \sum_{i\in I}\left| X_{i,n}\right|$. As  $Q$ belongs to the dual of $L^1(c)$, $E_Q\left[  \sum_{i\in I}\left| X_{i,n}\right|\right]<\infty$ and it follows from the dominated convergence theorem, that 
$$\lim_{k\to \infty}E_Q\left[ -\sum_{i\in I}X_{i,n}f_{i,n}^k \right]=E_Q\left[ -Z_n \right]$$
So, there exists $k\in\N$ such that,
$$ E_Q\left[\alpha_{s,t}^m\left( Q\right) \right] \leqslant E_Q\left[-Z_n \right] +\epsilon\leqslant E_Q\left[ -\sum_{i\in I}X_{i,n}f_{i,n}^k \right]+2\epsilon$$

Thus, using the fact that $\sum_{i\in I}X_{i,n} f_{i,n}^k\in \mathcal{A}_{s,t}$,
$$ E_Q\left[E_Q\left[\alpha_{s,t}^m\left( Q\right)\vert \mathcal{B}_r \right]\right] \leqslant E_Q\left[\makebox{ess sup}_{Z\in\mathcal{A}_{s,t}}E_Q\left[ -Z\vert \mathcal{B}_r\right]\right]$$
And the  equality (\ref{eqE}) is proved.
\item
If $E_Q\left[\alpha_{s,t}^m\left( Q\right)\right]=+\infty$, $\left( E_Q\left[-Z_n \right] \right)_{n\in\N}$ is increasing and tends to infinity. So, using the same ideas as before, we get the result. \hfill $\square $
\end{itemize}

\subsection{Regularity of paths}
Notice that for every $s \leq T$, the knowledge of  $\rho_{sT}(X)$ $P$ a.s. for some probability measure $P$ belonging to the canonical class determines entirely the element $\rho_{sT}(X)$ of $L^1_s(c)$.  Therefore it is very important to prove the existence of a process $Y_s$ with regular paths such that for every $s$, $ \rho_{sT}(X)=Y_s \;P$ a.s. This is the object of the following proposition

\begin{theorem}
Let $(\rho_{s,t})$ be a normalized time consistent convex dynamic risk measure on $L^1(c)$. Let $P$ be a probability measure belonging to the canonical $c$-class. Assume  that  $\alpha^m_{0T}(P)=0$.  Then for every $X \in L^{1}_T(c)$, $(\rho_{s,T}(X))_s$ is a $P$-supermartingale. \\Furthermore for every bounded $X$  in $L^1(c)$, the map $s \rightarrow E_P(\rho_{s,T}(X))$ is rightcontinuous. For every  $X$  bounded  $(\rho_{s,T}(X))_s$ admits $P$ a.s. a right continuous version which is furthermore optional. For every  $X$  bounded from above $(\rho_{s,T}(X))_s$ admits $P$ a.s. an optional version.
\label{propreg2}
\end{theorem}
Notice that even if the filtration is not complete we obtain an optional process.\\
{\bf Proof}  Let $0 \leq s \leq s' \leq T$. From the cocycle condition and the non negativity of the minimal penalty which are satisfied for every normalized time consistent dynamic risk mesure,  $\alpha^m_{ss'}(P)=0$. It follows from the time consistency that for every $X \in L^1(c)$, $\rho_{sT}(X) \geq E_P(\rho_{s'T}(X)|{\cal B}_s)$, i.e. $(\rho_{s,T}(X))_s$ is a $P$-supermartingale.\\
The most tricky part of the proof consists in proving the right continuity of the map $s \rightarrow E_P(\rho_{s,T}(X))$ for $X$ bounded.  To prove the right continuity we will adapt the proof  which   was  given in \cite{JBN} for dynamic risk measures on $L^{\infty}$ spaces under the hypothesis of time consistency for stopping times. Therefore  we explain here which changes have to be done. 
Let $s<t$. From the proof of $iii) \Rightarrow i)$ of Theorem \ref{thmtc}, the set 
$\{E_{R}(-X|{\cal B}_s)- \alpha^m_{st}(R),\; R \in {\cal Q}\}$ is  a lattice upward directed.  Thus from \cite{nev}, there is a sequence $\{R_n,\; n \in \N\}$ of probability measures in the dual of $L^1(c)$ whose restriction to ${\cal B}_s$ is equal to $P$   such that 
$\rho_{sT}(X)$ is the increasing limit of $(E_{R_n}(-X|{\cal B}_s)- \alpha^m_{sT}(R_n))$. Thus for every $k \in \N^*$, there is $S_k$ among the $R_n$ such that 
$E_P(\rho_{sT}(X))-\frac{1}{k} \leq (E_{S_k}(-X)- E_P(\alpha^m_{sT}(S_k))$. The proof follows then the first step of the proof of Lemma 4 of \cite{JBN}, in particular making use of the cocycle condition for $S_k$, and of the inequality $
(E_{S_k}(-X|{\cal B}_{s_n})- \alpha^m_{s_nT}(S_k)) \leq \rho_{s_nT}(X)\;\;P \;a.s.$, 
we get 
for a sequence $s_n<T$ decreasing to $s$,
\begin{equation}
E_P(\rho_{sT}(X))-\frac{1}{k} \leq E_{S_k}(\rho_{s_nT}(X))- E_P(\alpha^m_{ss_n}(S_k))
\label{eqcad1}
\end{equation}
Notice that the main difference with the proof of Lemma 4 of \cite{JBN} is that in general $S_k$ is not absolutely continuous with respect to $P$ therefore  we introduce a new probability measure $\tilde S_k=\frac{k-1}{k}P+\frac{1}{k}S_k$. It follows from (\ref{eqcad1}) that 
\begin{eqnarray}
E_P(\rho_{sT}(X))-\frac{1}{k} \leq & E_{P}(\rho_{s_nT}(X))- E_P(\alpha^m_{ss_n}(S_k))\nonumber \\ & +E_{\tilde S_k}(\rho_{s_nT}(X))[E_{\tilde S_k}(\frac{dS_k}{d\tilde S_k}|{\cal B}_{s_n})-E_{\tilde S_k}(\frac{dP}{d\tilde S_k}|{\cal B}_{s_n})]
\label{eqcad2}
\end{eqnarray}
The restrictions of $S_k$ and $\tilde S_k$  to ${\cal B}_s$ are both equal to $P$, then  $E_{\tilde S_k}(\frac{dS_k}{d\tilde S_k}|{\cal B}_{s_n})-E_{\tilde S_k}(\frac{dP}{d\tilde S_k}|{\cal B}_{s_n}) \rightarrow 0$ in $L^1(\tilde S_k)$ as $n \rightarrow \infty$.\\
As $(\rho_{s_n,T}(X))_{n \in \N}$ are uniformly bounded, it follows from the dominated convergence theorem that $E_{\tilde S_k}(\rho_{s_nT}(X))[E_{\tilde S_k}(\frac{dS_k}{d\tilde S_k}|{\cal B}_{s_n})-E_{\tilde S_k}(\frac{dP}{d\tilde S_k}|{\cal B}_{s_n})] \rightarrow 0$ as $n \rightarrow \infty$. Thus from equation (\ref{eqcad2}), there is $N(k)$ such that for $n \geq N(k)$, $E_P(\rho_{sT}(X))-\frac{1}{k} \leq  E_{P}(\rho_{s_nT}(X))- E_P(\alpha^m_{ss_n}(S_k))+\frac{1}{k}$. Due to the normalization, $\alpha^m_{ss_n}(S_k) \geq 0$. 
As  $\rho_{sT}(X)$ is a $P$-supermartingale, it follows that  $E_P(\rho_{s_nT}(X)) \leq E_P(\rho_{sT}(X))$. \\
This proves that $E_P(\rho_{s,T}(X)$ is the limit of $E_P(\rho_{s_n,T}(X))$.  
Then the proposition  for $X$ bounded follows from the modification Theorem of \cite{DM} (Theorem 4 page 76), and more precisely from the remark 5 following this Theorem. Notice that   as the filtration ${\cal B}_s$ is right continuous but not complete, the right continuity does not imply the existence of an optional version, but the two results are proved in  \cite{DM}. When $X$ is bounded from above, $X$ is the limit of the decreasing sequence  $X_n=sup(X,-n)$. From the theorem of representation, for every $s$, $\rho_{st}(X)$ is $P$ a.s. the increasing limit of $\rho_{st}(X_n)$. It is thus optional.   \hfill $\square $\\
From the preceding theorem we  deduce easily the following result. 
 \begin{corollary}
 Let $\rho_{s,t}$ be a normalized time consistent convex dynamic risk measure on $L^1(c)$. Assume that for every probability measure $Q$ in $K_+$  (the non negative part of the unit ball of the dual of $L^1(c)$), there is a probability measure $P$ belonging to the canonical $c$-class  such that $Q \ll P$ and $\alpha^m_{0T}(P)=0$.\\
Then for every $X \in L^{1}_T(c)$, and $Q \in K_+$, $(\rho_{s,T}(X))_s$ is a $Q$-supermartingale. 
For every bounded $X$ for every choice $Y_s$ of $\rho_{sT}(X)$, there is a right continuous process $Z_s$ such that for every $Q \in K_+$, 
 $Q(\{Y_s \neq Z_s\})=0$. 
 \label{correg1}
 \end{corollary}
 {\bf Proof}
 let $0 \leq s \leq T$.  let $Y_s$ be any choice in the $c$-class of $\rho_{sT}(X)$. Let $A=\{\omega\;|\;\exists s \in [0,T[,
Y_s(\omega) \neq \limsup_{\substack{t>s\\t \rightarrow s}} Y_t(\omega) \;or Y_s(\omega)
\neq \liminf_{\substack{t>s\\t \rightarrow s}} Y_t(\omega)\}$. From the above proposition, for every $P$ in the canonical $c$-class such that $\alpha^m_{0T}(P)=0$, $P(A)=0$
Thus $Q(A)=0$ for every $Q \in K_+$. The process $ (Z_s)_{0 \leq s \leq T}$ defined as $Z_s=1_{A^c}Y_s$ satisfies obviously the required conditions. \hfill $\square $
 
The preceding result applied to the  particular case of sublinear risk measures gives:
 \begin{corollary}
 Let $(\rho_{s,t})_{0 \leq s \leq t \leq \infty}$ be a sublinear time consistent dynamic risk measure on $L^1(c)$. Let ${\cal Q}$ be a weakly relatively compact set  of probability measures such that $\rho_{0,\infty}(X)=\sup_{Q \in{\cal Q}}E_Q(-X)$. Then for every $X$ in $L^1(c)$, $\rho_{sT}(X)$ is a $Q$-supermartingale for every $Q \in {\cal Q}$. Furthermore for every $X$ bounded from above and  every choice $Y_s$ of $\rho_{sT}(X)$, there is a right continuous process $Z_s$ such that for every $Q \in {\cal Q}$, 
 $Q(\{Y_s \neq Z_s\})=0$. 
\label{correg2}
 \end{corollary}
{\bf proof}
 The existence  of a weakly relatively compact set of probability measures such that $\rho_{0,\infty}(X)=\sup_{Q \in{\cal Q}}E_Q(-X)$ follows from \cite{BNK}, Proposition 3.1. We apply then the result of corollary \ref{correg1} with the capacity $c'(X)=\sup_{Q \in{\cal Q}}E_Q(|X|)$ \hfill $\square $\\
 Notice that  a  statement similar to that of  Corollary \ref{correg2} is enounced  in \cite{NS} for the filtration  completed with all polar sets. 
 
 \section{Sublinear Risk measures on ${\cal C}_b(\Omega)$ factorizing through time}

\subsection{Factorization through time} 
In this section $\Omega$ is a general Polish space.
 In \cite{BNK}, section 5,  we have introduced the notion of  regularity for a sublinear  risk measure $\rho$ on ${\cal C}_b(\Omega)$ which means that for every decreasing sequence $f_n$ with limit $0$, the sequence $\rho(-f_n)$ tends to  $0$. We have proved  in \cite{BNK} that  every regular sublinear risk measure on ${\cal C}_b(\Omega)$ is represented by a weakly relatively compact set of probability measures ${\cal P}$, i.e.  
 \begin{equation}
\rho(X)=\sup_{P \in {\cal P}}E_P(-X)\;\; \forall X \in {\cal C}_b(\Omega)
\label{eqB1}
\end{equation}
We address now the following question: 
 Can we find a sufficient condition such that the risk measure $\rho$ can be factorized through time, i.e. such that there is a time consistent dynamic risk measure $\rho_{st}$ satisfying $\rho_{0,\infty}=\rho$?  
In the case of risk measures on $L^{\infty}$ spaces,   Delbaen \cite{D} has introduced the 
notion of m-stability of a set of probability measures in order to answer this question. Notice that for every sublinear time consistent dynamic risk measure $\rho_{st}$, for every $s \leq t$, $\rho_{st}$ is equal to the restriction of $\rho_{s \infty}$, we can thus drop the second index, and a time consistent sublinear risk measure is thus a family $\rho_s$ such that $\rho_s=\rho_s(-\rho_t)$.\\
 Let $({\cal F}_t)$ be a right continuous filtration such that  ${\cal F}_{\infty}= B(\Omega)$ and ${\cal F}_0$ is trivial (i.e. is contained 
  in  the  $P$-null sets).  It is natural to address the question of the existence of a factorization on $L^{\infty}(\Omega,{\cal B}(\Omega), ({\cal F}_t), P)$ where $P$ is a probability measure belonging to  the canonical $c_{\rho}$ class ($c_{\rho}(X)=\rho(-|X|)$).  In the dual representation (\ref{eqB1}) of the risk measure $\rho$, one can replace  the set ${\cal P}$ by the weak closure of its convex hull. The sublinear risk measure is unchanged. Thus we define now the notion of stability for a convex weakly closed set of probability measures.
\begin{definition}
Let ${\cal C}$ be a  
convex weakly closed set of possibly non dominated probability measures on $(\Omega, {\cal B}(\Omega))$. 
${\cal C}$ is stable if for all probability measures $P,Q$ in ${\cal C}$, such that $Q \ll  P$, for every stopping time $\tau$, there is a probability measure $R \in {\cal C}$, $R \ll P$, such that $(\frac{dR}{dP})_t= (\frac{dQ}{dP})_t\;\forall t \leq \tau$  and  $(\frac{dR}{dP})_t=(\frac{dQ}{dP})_{\tau} \;\forall t \geq \tau$.
\label{defstable}
\end{definition}

\begin{proposition}
Let ${\cal C}$ be a  stable 
convex weakly compact set of possibly non dominated probability measures on ${\Omega}$.
Let $\rho(X)=\sup_{Q \in {\cal C}} E_Q(-X) \;\; \forall X \in {\cal C}_b(\Omega)$.\\ Then  there is a probability measure $P \in {\cal C}$ belonging to the canonical $c_{\rho}$-class, and a set ${\cal S}$ of probability measures in $L^1(P)$, m-stable with respect to the right continuous filtration ${\cal F}_{t}$, such that 
  $\rho(X)=\sup_{Q \in {\cal S}} E_Q(-X)$.\\
The dynamic risk measure $\rho_{\sigma}$ defined on $L^{\infty}(\Omega, {\cal B}(\Omega),P)$ by
$$\rho_{\sigma}(X)=\esssup_{\{Q \in {\cal S}\;Q \sim P\}} E_Q(-X|{\cal F}_{\sigma})\;\;P\;a.s.$$
is a time consistent dynamic risk measure on  $L^{\infty}(\Omega, {\cal B}(\Omega),({\cal F}_t)_{t \in \R_+},P)$ which factorizes  $\rho$, i.e. for every $X \in {\cal C}_b(\Omega)$,  $\rho(X)=\rho_{0}(X)$. 
\label{thmstable}
\end{proposition}
{\bf Proof}
From \cite{BNK},  there is a countable set $\{Q_n, \; n \in \N \}$ of probability measures in ${\cal C}$ such that 
 $\rho(X)=\sup_{n \in \N} E_{Q_n}(-X) \;\; \forall X \in L^1(c)$.
Let $P= \sum_{n \in \N} \frac{Q_n}{2^{n+1}}$. From \cite{BNK}, $P$ belongs to the canonical $c_{\rho}$-class.  As  ${\cal C}$ is convex and weakly compact, $P$ belongs to ${\cal C}$. Let ${\cal S}={\cal C} \INTER \{Q \ll  P\}$. ${\cal S} \subset L^1(P)$. The stability of ${\cal C}$ implies the m-stability of ${\cal S}$ for the right-continuous filtration $({\cal F}_{t})$ (cf \cite{D} for the definition of m-stability). For all 
$X$ in ${\cal C}_b(\Omega)$, 
$$\sup_{Q \in {\cal S}}E_Q(-X) \leq \sup_{Q \in {\cal C}}E_Q(-X)= \sup_{n  \in \N}E_{Q_n}(-X)\leq \sup_{Q \in {\cal S}} E_Q(-X)$$
Thus 
\begin{equation}
\forall X \in {\cal C}_b(\Omega), \;\rho(X)=\sup_{Q \in {\cal S}}E_Q(-X)
\label{eqext0}
\end{equation}
Furthermore every $Q$ in ${\cal S}$ is the weak limit of $Q_n=
(1-\frac{1}{n})Q+\frac{1}{n}P$,     $Q_n \in {\cal S}$, and $Q_n \sim P$. 
Thus 
\begin{equation}
\forall X \in {\cal C}_b(\Omega), \;\rho(X)=\sup_{Q \in {\cal S}\;Q \sim P}E_Q(-X)
\label{eqext}
\end{equation}
From \cite{D} the dynamic risk measure  $ \rho_{\sigma}$ defined on $L^{\infty}(\Omega, {\cal B}(\Omega),({\cal F}_t)_{t \in \R_+}, P)$ by  $\rho_{\sigma}(X)=\esssup_{\{Q \in S, \;\;Q \sim P \}}E_Q(-X|{\cal F}_{\sigma})\;\;P \;a.s.$ is time consistent. From (\ref{eqext}) it    extends the risk measure $\rho$ on ${\cal C}_b(\Omega)$. \hfill $\square $\\
It follows from the proof  that the conclusion of Proposition \ref{thmstable} is satisfied for every probability measure $P=\sum_{n \in \N}\alpha_n Q_n$, where $\{Q_n\}$ is a dense subset of ${\cal C}$.
\begin{corollary}
Let ${\cal C}$ and $\rho$ be as in Proposition 
\ref{thmstable}. Let $c_{\rho}$ be the capacity associated to $\rho$. The risk measure $\rho$ extends uniquely to $L^1(c_{\rho})$ and $\rho(X)=\sup_{Q\in {\cal S}} E_Q(-X)$ for every $X \in L^1(c_{\rho})$. For every stopping time $\sigma$, $\rho_{\sigma}$ has a unique continuous extension to $L^1(c_{\rho})$.
\label{corsta}
\end{corollary}
{\bf proof}
For all f in ${\cal C}_b(\Omega)$, $c_{\rho}(f)= \rho(-|f|)$.
 The risk measure $\rho$ extends uniquely to $L^1(c_{\rho})$.
From  \cite{BNK} the equalities $\rho(X)=\sup_{Q\in {\cal S}} E_Q(-X)$ and $c_{\rho}(X)=\sup_{Q \in {\cal S}} E_Q(|X|)$ are  also satisfied for every $X \in L^1(c_{\rho})$.
From the m-stability of ${\cal S}$, it follows that for every $f \in {\cal C}_b(\Omega)$, the set $\{E_Q(|f||{\cal F}_{\sigma}),\; Q \in {\cal S},\;\; Q \sim P\}$ is  a lattice upward directed. Thus   $\sup_{R \in {\cal S}}E_R(|\rho_{\sigma}(f)|) \leq \sup_{S \in {\cal S}}E_S(|f|)=c_{\rho}(f)$. It follows that  for every stopping time $\sigma$, $\rho_{\sigma}$ has a unique continuous extension to $L^1(c_{\rho})$ with values in $\cap_{Q \in {\cal S}}L^1(\Omega,{\cal F}_t,Q)$.
\begin{remark}
Even if ${\cal F}_t$ is the right continuous filtration generated by the open sets of $\Omega_t$, it can happen that  $\rho_{t} (X)$ does not belong to $L^1_t(c_{\rho})$, and that the time consistency property which is satisfied for all continuous bounded functions cannot be extended to $L^1(c_{\rho})$.
\label{rq1}
\end{remark}
Therefore we introduce now a more restrictive notion of factorization through time for a sublinear risk measure defined on ${\cal C}_b(\Omega)$. 

\subsection{Factorization on $L^1(c_{\rho})$}
 We restrict now to the case where $\Omega$ is either equal to ${\cal C}_0 ([0,\infty[,\R^d)$ or $D([0,\infty[,\R^d)$, and ${\cal F}_t={\cal B}_t$ is the right continuous filtration generated by the open sets of $\Omega_t$. It seems to be too restrictive to ask for the existence of a  factorization  within the set of continuous bounded functions. Thus we study now factorization on $L^1(c_{\rho})$.

\begin{definition}
Let $\rho$ be a sublinear regular   risk measure  on ${\cal C}_b(\Omega)$. Let $c_{\rho}$ be the associated capacity: $c_{\rho}(X)=\rho(-|X|)$. We say that $\rho$ factorizes through time on $L^1(c_{\rho})$, if there is a time consistent sublinear dynamic risk measure $\rho_{st}=\rho_s$ on $L^1(c_{\rho})$ such that 
\begin{equation}
\forall f \geq 0, f  \in {\cal C}_b(\Omega_s), \forall X \in L^1(c_{\rho}), 
\rho_{s}(fX) \leq f \rho_s(X)
\label{eqss}
\end{equation}
extending $\rho$ (i.e. $\forall X \in  {\cal C}_b(\Omega)$, $\rho(X)=\rho_{0}(X)$). 
\label{deffac}
\end{definition}
By definition of a dynamic risk measure on $L^1(c_{\rho})$, necessarily for every  $X$ in $L^1(c_{\rho})$, $\rho_s(X) $ belongs to $L^1_s(c_{\rho})$. 
\begin{proposition}
Let $\rho$ be a sublinear regular risk measure on ${\cal C}_b(\Omega)$ factorizing  through time on $L^1(c_{\rho})$ into $\rho_{s}$. Let $P$ in $K_+$ (the non negative part of the unit ball of the dual of $L^1(c_{\rho})$) belonging to the canonical $c_{\rho}$-class. Then for every $s \geq 0$,
\begin{itemize}
\item There is a countable  set $\{Q_n, n \in \N\}$ of probability measures belonging to  the unit ball of the dual of $L^1(c_{\rho})$, whose restriction to ${\cal B}_s$ is equal to $P$ such that 
\begin{equation}
\forall X\in L^1(c_{\rho}), \rho_{s}(X)=P \esssup_{ n \in \N} E_{Q_n}(-X|{\cal B}_s)\; P\;a.s.
\label{eqrep20b}
\end{equation}
\item For every $X$ in $L^1(c_{\rho})$, there is a probability measure $Q_X$ belonging to the unit ball of the dual of $L^1(c_{\rho})$ whose restriction to ${\cal B}_s$ is equal to $P$ such that  
\begin{equation}
 \rho_{s}(X)= E_{Q_X}(-X|{\cal B}_s)\;\; P\;a.s.
\label{eqrep21b}
\end{equation}
\end{itemize}
\label{propfac}
\end{proposition}
{\bf Proof}
For every $s$,  $\rho_s$ satisfies the strong convexity property, thus from Theorem \ref{thm3} there is a sequence $Q_n$ of probability measures on $(\Omega,{\cal B}(\Omega))$, belonging to the dual of $L^1(c_{\rho})$ whose restriction to ${\cal B}_s$ is equal to $P$ satisfying  (\ref{eqrep20b}).
From the monotonicity and time consistency condition, $c_{\rho}(\rho_{s}(-|X|))=\rho(-|X|)=c_{\rho}(|X|)$ . Thus for every $n$ and $X \in L^1(c_{\rho})$,\\
 $E_{Q_n}(|X|)=E_P(E_{Q_n}(|X||{\cal B}_s) \leq E_P(\rho_{s}(-|X|))) \leq c_{\rho}(|X|)$. This proves that every $Q_n$ belongs to the unit ball of  the dual of $L^1(c_{\rho})$.\\
The second part of the proposition  follows from corollary \ref{c1} and the above argument applied with $Q_X$.  \hfill $\square $

We give now a sufficient condition for a sublinear risk measure on 
${\cal C}_b(\Omega)$ to admit a factorization through time on $L^1(c_{\rho})$.

\begin{theorem}
Let ${\cal C}$ be a stable convex weakly compact set of possibly non dominated probability measures on $\Omega$ such that $\rho(f)=\sup_{Q \in {\cal C}}E_Q(-f)$. Let $P$ in $K_+$ belonging to the canonical $c_{\rho}$-class.  Assume  that for $X$ in a dense subset of ${\cal C}_b(\Omega)$ (for the $c$-norm),  there is an element $\rho_{s}(X)$ of $L^1_s(c_{\rho})$ such that 
\begin{equation}
\rho_{s}(X)=\esssup_{\{Q \in {\cal S}\; Q \sim P\}} E_Q(-X|{\cal F}_{s}) \;\;P\;a.s.
\label{eqpro}
\end{equation}
Then for all $s$, $\rho_{s}$ extends to  $L^1(c_{\rho})$ with values in $L^1_s(c_{\rho})$.
 The sublinear risk measure $\rho$   on 
${\cal C}_b(\Omega)$ factorizes through time on $L^1(c_{\rho})$ (cf definition (\ref{deffac})). In particular 
\begin{equation}
 \forall r \leq s, \;\forall X \in L^1(c_{\rho}),\; \rho_{r}(X)=\rho_{r}(-\rho_{s}(X))\;\; 
\label{eqstc}
\end{equation} 
\label{thmfac}
\end{theorem}
The above    equation is an equality in $L^1(c_{\rho})$. 

{\bf Proof}
From Proposition \ref{propfac}, equation (\ref{eqpro}) defines a time consistent dynamic risk measure $\rho_s$ on $L^{\infty} (\Omega,{\cal B}(\Omega),P)$ factorizing $\rho$.  
It follows easily from sublinearity and monotonicity of $\rho_{s}$ that for every $X,Y \in {\cal C}_b(\Omega)$,  $|\rho_{s}(X)-\rho_{s}(Y)|\leq \rho_{s}(-|X-Y|)$. From Lemma 4.1 of \cite{BNK}, it follows that on the non positive  elements of $L^1(c_{\rho})$, $c_{\rho}$ is equal to $\rho_0 $. Denote ${\cal D}_s$ a dense subset of ${\cal C}_b(\Omega)$ such that for all $X \in {\cal D}_s$, $\rho_s(X) \in L^1_s(c_{\rho})$. From the monotonicity of $c_{\rho}$ and the time consitency it  follows that for all $X,Y \in{\cal D}_s$,
\begin{equation}
c_{\rho}(|\rho_{s}(X)-\rho_{s}(Y)|)= \rho_0(-|\rho_{s}(X)-\rho_{s}(Y)|)
 \leq \rho_{0}(-\rho_s(-|X-Y|) =c(|X-Y|)
\label{eqexte}
\end{equation}
 This proves that $\rho_{s}$ extends uniquely to $L^1(c_{\rho})$ with  values in $L^1_s(c_{\rho})$. 
For all $r \leq s \leq t$ and every $X \in {\cal C}_b(\Omega)$, $\rho_{r}(X)$ and $\rho_{r}(-\rho_{s}(X))$ belong to $L^1(c_{\rho}))$. From Proposition \ref{thmstable} they  are equal $P$ a.s. As $P$ belongs to the canonical $c_{\rho}$ class, this means that $\forall X \in{\cal C}_b(\Omega)$,
\begin{equation}
\;\;\rho_{r}(X)=\rho_{r}(-\rho_{s}(X))
\label{eqext2}
\end{equation}
 The equality (\ref{eqext2}) is an equality in $L^1(c_{\rho})$. From the equation (\ref{eqexte}), it then follows that the equality (\ref{eqext2}) is satisfied for all $X \in L^1(c_{\rho})$.  \hfill $\square $
 \begin{remark}
 \begin{itemize}
 \item
Nutz and Soner address in \cite{NS} the question of unicity of a time consistent extension of a sublinear risk measure, given a set of probability measures and given a family of vector spaces ${\cal H}_t \subset {\cal H}$.
It follows from Proposition 3.1 of \cite{NS} that both factorization through time on $L^{\infty}(\Omega,{\cal B}(\Omega),{\cal F}_t,P)$ where $P$ belongs to the canonical $c_{\rho}$-class and factorization on $L^1(c_{\rho})$ are unique.
\item
Notice that in all the existing examples of time consistent dynamic risk measures in case of uncertainty represented by a weakly relatively compact set of probability measures, the dynamic risk measure $\rho_t$ is defined on a vector space ${\cal H}$ which is always a subspace of $L^1(c_{\rho})$. There is never unicity of a set of probability measures ${\cal P}$ such that $\rho(X)= \sup_{P \in {\cal P}} E_P(-X)\;\; \forall X  \in {\cal C}_b(\Omega)$. As proved in \cite{BNK} there is always a countable set ${\cal P}$ such that the above equality is satisfied. On the other hand  the set of all probability measures $Q$ such that $E_Q(-X) \leq \rho(X)$ for all $X$   in $L^1(c_{\rho})$ is convex and weakly compact. The set of probability measures which is choosen in all the papers constructing examples of time consistent dynamic risk measures \cite{STZ} and \cite{NS} is an  intermediary set: not countable nor convex.
\item
Notice that the factorization through time on $L^1(c_{\rho})$ is in some sense universal. Indeed  when it exists (cf a sufficient condition for existence in Theorem \ref{thmfac}), it does not depend on the choice of a particular set of probability measures representing $\rho$. The Banach spaces $L^1_t(c_{\rho})$ depend only on $c_{\rho}(f)=\rho(-|f|)$. The equality $\rho_s=\rho_s \circ(-\rho_t)$ is satisfied in $L^1(c_{\rho})$ and thus it is satisfied $Q\; a.s.$ for every $Q \in {\cal P}$ whatever  the choice of the set ${\cal P}$ representing $\rho$.  
\end{itemize}
\end{remark}
\section{Examples of dynamic risk measures}
In this section we restrict to the particular case  where  $\Omega={\cal C}_0([0,T],\R^{d})$ the space of continuous functions on $[0,T]$ null in zero. Let ${\cal B}(\Omega)$ be the Borel $\sigma$-algebra. let $\Omega_t= {\cal C}_0([0,t],\R^{d})$.  Let ${\cal B}_t$ be the right continuous filtration  generated by the open sets of $\Omega_t$. Denote $B_t$ the coordinate process. 
\subsection{Dynamic pricing in context of uncertain volatility}
 \label{secvolu}
 We consider the framework  introduced in \cite{DM}. A probability measure $Q$ on $(\Omega, {\cal B}(\Omega),{\cal B}_t)$ is called a martingale measure if the coordinate process $(B_t)$ is a martingale with respect to ${\cal B}_t$ under $Q$ and if the martingales
$(( B_{i})_t)_{1\leq i\leq d}$ are orthogonal in
the sense that for all $i\neq j$, $ <B_{i},B_{j}>^Q_t=0\ \ Q\ a.s. $ where  $<B_{i}, B_j>^Q$ denotes the quadratic covariational process
corresponding to  $B^{i}$ and $B^j$, under $Q$ and $<B>^Q$ the quadratic variation of $B$ under $Q$. Fix for all $i\in\left\{ 1,\dots,d\right\}$ two finite deterministic H\"{o}lder-continuous measures $\underline{\mu}_{i}$ and $\mu_{i}$ on $[0,T]$  and consider the set ${\cal P}$ of orthogonal martingale measures such that $$\forall i\in\{1,\dots,d\},\ \ \   d\underline{\mu}_{i,t} \leq d<B_{i}>_{t}^{Q}
\leq d\mu_{i,t}.$$ From  \cite{BNK} the set ${\cal P}$ is convex and weakly compact. Let $\rho$ be the risk measure  defined on ${\cal C}_b(\Omega)$ by $\rho(f)= \sup_{P \in {\cal P}} E_P(-f)$, and $c_{\rho}$ the associated capacity, $c_{\rho}(f)=\rho(-|f|)$. 
\begin{proposition}
 There is a probability measure $P$ belonging to the canonical $c_{\rho}$-class and a m-stable  set $S$  of probability measures all absolutely continuous with respect to $P$ such that the dynamic risk measure  on $L^{\infty}(\Omega,{\cal B}_t,P)$ defined by    $\rho_{s} (X)=\esssup_{P \in S} E_P(-X|{\cal B}_s)\;P \;a.s.$  is time consistent and factorizes $\rho$ (i.e.$\rho_{0}=\rho$).
\label{propDM}
\end{proposition}
{\bf Proof}
 From  \cite{BNK} the set ${\cal P}$ is convex and weakly compact.  \\
From \cite{DM}, $(B_i)_s^2$ belongs to  $L^1(c_{\rho})$ for every $s$, and thus  the quadratic variation of $B$ is defined as an element of  $L^1(c_{\rho})$ by the equation
\begin{equation}
<B_i>^{c}_t=(B_i)_t^2-2 \int_0^t (B_i)_s d(B_i)_s
\label{qv}
\end{equation} In the same way, $<B_i,B_j>$ is defined in $L^1(c_{\rho})$. 
Let $S$ and $Q$ be two  probability measures in ${\cal P}$ such that $Q \ll S$. The probability measure $R$ such that $(\frac{dR}{dS})_t=(\frac{dQ}{dS})_t \;\forall t \leq \tau$  and  $(\frac{dR}{dS})_t=(\frac{dQ}{dS})_{\tau} \;\forall t \geq \tau$ is also absolutely continuous with respect  to $S$. For all $i$,  the inequality  $d\underline{\mu}_{i,t} \leq d<B_{i}>_{t}^{S}
\leq d\mu_{i,t}$ is satisfied  $S\; a.s.$ , thus is also satisfied $R$ a.s.. In the same way, $<B_i,B_j>=0\;R\; a.s.$. On the other hand,  the set of  martingale measures for one process (here $B_t$) is m-stable (cf \cite{D}), thus $B_t$ is a $R$ martingale. Thus $R$ belongs to ${\cal P}$ and this proves that ${\cal P}$ is stable. The result follows then from Proposition \ref{thmstable}.  \hfill $\square $

From Proposition \ref{propDM} and Remark \ref{rq1}, we deduce the following result:
\begin{corollary}
Let ${\cal P}$ be as above. To every  $X$ in ${\cal C}_b(\Omega_t)$ one can associate its  dynamic ask price 
$$\Pi_{st}(X)=\esssup_{P \in {\cal S}}E_P (X|{\cal B}_s)$$
and its dynamic bid price 
$$-\Pi_{st}(-X)=\essinf_{P \in {\cal S}}E_P (X|{\cal B}_s)$$ The above formulas  define a time consitent dynamic pricing procedure on $L^{\infty}(\Omega,({\cal B}_t)_{t \in \R^+},P)$ in the sense of \cite{JBN4}, with no arbitrage. The above formulas extend to $X \in L^1(c)$ ($c(X)=\sup_{P \in {\cal P}}E_P (|X|)$).  Furthermore for every $X \in L^1(c)$,  $\Pi_{0T}(X)=\sup_{P \in {\cal P}}E_P (X)$.
\label{corTCPP}
\end{corollary}
\subsection{Conditional G-Expectations}
 \label{G}
Notice that in the  preceding section we have constructed  an example of a sublinear risk measure on ${\cal C}_b(\Omega)$  admitting an extension into a time consistent dynamic risk measure $\rho_s$ on $L^{\infty}(\Omega,{\cal B}_s,P)$. However in the previous example it can happen that  for some $X \in {\cal C}_b(\Omega)$, $\rho_{s}(X)$ does not belong to $L^1_s(c)$.  We prove now that the conditional  $G$-expectations takes always values in $L^1_s(c)$ and defines thus a sublinear risk measure on ${\cal C}_b(\Omega)$ factorizing through time on $L^1(c)$ (cf definition(\ref{deffac})).\\    
In  this section, $\Omega={\cal C}_0 ([0,\infty[,\R^d)$.     
Peng introduced the notion  of  G-expectations (\cite{P1}and  \cite{P2}) defined on $Lip=\cup_T Lip_T$.\\  $Lip_T=\{\phi(B_{t_1},B_{t_2},...B_{t_k}),\; t_1 \leq t_2 \leq t_k \leq T|\;\phi \in {\cal L}ip(\R^k)\},\;k \in \N_+^*\}$ where ${\cal L}ip(\R^k)$ denotes the set of $\R$ valued Lipschitzian functions on $\R^k$. \\ 
G-expectations are defined from solutions of P.D.E.  For a complete study on G-expectations  we refer to \cite{P3}. Here we will use the notations of \cite{P3}. For $0 \leq t$, $\hat E (.|{\cal F}_t):Lip \rightarrow Lip_t$ is the  conditional G-expectation. For every $X \in Lip$, let $\rho_{s}(X)=\hat E (-X|{\cal F}_s)$ and $c(X)=\hat E (|X|)$
\begin{proposition}
$c(X)=\hat E (|X|)$ defines a capacity  on ${\cal C}_b(\Omega)$.
For every $0 \leq s \leq t $, $\hat E (.|{\cal F}_s)$ extends uniquely to $L^1(c)$ with values in $L^1_s(c)$. This extension still denoted $\hat E (.|{\cal F}_s)$ defines (up to a minus sign) a sublinear strongly convex time consistent dynamic risk measure on $L^1(c)$.
In particular there is a countable set $\{Q_n, n \in \N\}$ of probability measures  whose restriction to ${\cal F}_s$ is equal to $P$ such that  $\hat E (.|{\cal F}_s)$ admits the following dual representation.
\begin{equation}
\forall X\in L^1(c),\;\hat E (X|{\cal F}_s)=P \esssup_{ n \in \N} E_{Q_n}(X|{\cal F}_s)\;\; P\;a.s.
\label{eqrep20ex}
\end{equation}
Furthermore for every $X \in L^1(c)$, there is a probability measure $Q_X$ in $K_+$  whose restriction to ${\cal F}_s$ is equal to $P$  such that 
\begin{equation}
\hat E (X|{\cal F}_s)=E_{Q_X}(X|{\cal F}_s)
 \label{eqfin}
 \end{equation}
 \end{proposition}
 
{\bf Proof}
The existence of a weakly relatively compact set of probability measures ${\cal P}$ such that  for every $f \in Lip$, 
$\hat E(f)=\sup_{P \in {\cal P}}E_P(f)\;\;$  is proved in \cite{DHP}. Thus $c(f)=\sup_{P \in {\cal P}}E_P(|f|)$ defines a capacity on ${\cal C}_b(\Omega)$. It is also proved in \cite{DHP} that $Lip$ is dense in ${\cal C}_b(\Omega)$ for the $c$-norm.
From the remark following Proposition 2.3 of Chapter III in \cite{P3} 
\begin{equation}
\forall  X,Y \in Lip, \;c(\hat E (X|{\cal F}_s)-\hat E (Y|{\cal F}_s))\leq c(X-Y)
\label{equnif}
\end{equation} 
As $\hat E (X|{\cal F}_s) \in Lip_s$ for every $X \in Lip$, it follows from equation (\ref{equnif}) that $\hat E (.|{\cal F}_s)$ admits a unique extension to $L^1(c)$ with values in $L^1_s(c)$. The monotonicity, sublinearity and time consitency  for $\hat E (.|{\cal F}_s)$ on $L^1(c)$ follow easily from proposition 2.3 of Chapter III in \cite{P3} using equation (\ref{equnif}).
From Proposition 2.3 of Chapter III in \cite{P3}, for every $f \geq 0$, $f \in Lip_s$ and $X \in Lip$,
$\hat E (fX|{\cal F}_s) =f \hat E (X|{\cal F}_s)$.  By density and continuity for the $c$ norm, it follows that this equality is satisfied for all  $f \geq 0$ in ${\cal C}_b(\Omega_s) $ and $X  \in L^1(c)$. (Notice that this equality is even satisfied  for every non negative bounded  $f$ in  $L^1_s(c)$).  The strong convexity follows then from this equality and the sublinearity. 
The last result follows from Corollary \ref{c1}. \hfill $\square $

\end{document}